\documentclass[11pt,a4paper]{amsart}
\usepackage{amsfonts,lscape}
\usepackage[top=33mm,right=33mm,bottom=33mm,left=33mm]{geometry}
\usepackage{amsmath,array}
\usepackage{amssymb,srcltx}

\numberwithin{equation}{section}
\newcommand{\GL}{\mathrm{GL}}
\newcommand{\SL}{\mathrm{SL}}
\newcommand{\SU}{\mathrm{SU}}
\newcommand{\GU}{\mathrm{GU}}
\newcommand{\Sp}{\mathrm{Sp}}

\newcommand{\SO}{\mathrm{SO}}
\newcommand{\PSO}{\mathrm{PSO}}
\newcommand{\POm}{\mathrm{P\Omega}}
\newcommand{\Om}{\mathrm{\Omega}}

\newcommand{\PGL}{\mathrm{PGL}}	

\newcommand{\PSp}{\mathrm{PSp}}
\newcommand{\Gbar}{\overline{G}}
\newcommand{\Tbar}{\overline{T}}
\newcommand{\GF}{G}

\newcommand{\tpo}{$2$-part order }
\newcommand{\tpe}{$2$-part exponent }
\newcommand{\tpes}{$2$-part exponents }
\newcommand{\tpos}{$2$-part orders }

\renewcommand{\l}{n}
\renewcommand{\a}{b}
\renewcommand{\j}{b}
\newcommand{\I}{\mathcal{I}}
\makeatletter
\def\imod#1{\allowbreak\mkern7mu({\operator@font mod}\,\,#1)}
\makeatother

\newtheorem{thm}{Theorem}[section]

\newtheorem{lemma}[thm]{Lemma}
\newtheorem{pro}[thm]{Proposition}
\newtheorem{cor}[thm]{Corollary}

\newtheorem{remark}[thm]{Remark}
\newtheorem{dfn}[thm]{Definition}
	
\title[Proportions of $2$-part order elements]{Proportions of elements with given $2$-part order in finite classical groups of odd characteristic}
\author{Simon Guest and Cheryl~E. Praeger}
\thanks{This research forms part of Australian Research Council Federation Fellow project FF0776186 of the second author.}
\address{Centre for the Mathematics of Symmetry and Computation (M019),
The University of Western Australia,
35 Stirling Highway,
Crawley WA 6009, Australia\footnote{Current address for S. Guest: School of Mathematics, University of Southampton, Southampton, SO17 1BJ, UK}}
\email{simon.guest@uwa.edu.au}
\email{cheryl.praeger@uwa.edu.au}

\date{\today}
\begin{document}
\maketitle
\begin{abstract}
For an element $g$ in a group $X$, we say that $g$ has \tpo $2^{\a}$ if $2^{\a}$ is the largest power of $2$ dividing the order of $g$.
We prove lower bounds on the proportion of elements in finite classical groups in odd characteristic that have certain \tpos \hspace{-1.2mm}.  In particular, we show that the proportion of odd order elements in the symplectic and orthogonal groups is at least $C/\l^{3/4}$, where $\l$ is the Lie rank, and $C$ is an explicit constant. We also prove positive constant lower bounds for the proportion of elements of certain \tpos independent of the Lie rank.
Furthermore, we describe how these results can be used to analyze part of Yal\c{c}inkaya's  Black Box recognition algorithm for finite classical groups in odd characteristic.
\end{abstract}

\section{Introduction}

The aim of this paper is to find lower bounds on the proportion of elements of certain \tpos in finite classical groups of odd characteristic.
Babai, Palfy, and Saxl  \cite{BPS} proved that for a prime $p$, the proportion of elements of order not divisible by $p$ in a finite classical \emph{simple} group $G$ is at least $1/2d$, where $d$ is the dimension of the natural module for $G$. In particular, the proportion of odd order elements in $G$ is at least $1/2d$. In this paper, we prove an improved lower bound if $p=2$ and $G$ is a symplectic or orthogonal group of odd characteristic, and also lower bounds for the proportions of elements of twice odd order in these groups.

 \begin{thm} \label{thm:oddandtwiceodd}
Let $q$ be an odd prime power, and suppose that $d \ge 4$. If $G=\Sp_{d}(q)$, $\SO^{\epsilon}_{d}(q)$, or $\Om_d^{\epsilon}(q)$ (where $\epsilon = \pm$ or $\circ$), then the proportion of elements in $G$ of odd order is at least $C_1/d^{3/4}$, and the proportion of elements in $G$ of twice odd order is at least $C_2/d^{1/2}$, where $C_1$ and $C_2$ are explicit constants. Moreover, these lower bounds also hold in $G/Z(G)$.
 \end{thm}

We prove Theorem \ref{thm:oddandtwiceodd}, and more general lower bounds for other \tpos and other classical groups, using recent work of Niemeyer and the second author \cite{NiePra}. This work exploits Lie theoretic methods first introduced by G.I. Lehrer \cite{Leh} to study characters of finite groups of Lie type. These methods have been used variously by Isaacs, Kantor and Spaltenstein \cite{IKS} to study the distribution of elements in permutation groups and by Cohen and Murray \cite{CohMur} to analyse algorithms for computing with simple Lie algebras. 
 We also make use of recent work of the authors \cite{GPsymmetric} on the proportion of elements of a given \tpo in a finite symmetric group $S_n$.

We denote by $k_2$ the $2$-part of a positive integer $k$; that is the largest power of $2$ dividing $k$.
Let  
\begin{equation} \label{e:a}
  a :=  \frac{1}{\log(3 + \sqrt{5}) - \log(2)} \approx 1.03904
\end{equation}
 and for an odd prime power $q$, let 
\begin{equation}\label{eqn:introdeft} 
  t(q) := \begin{cases}
  (q-1)_2, & \hbox{if $q \equiv 1 \imod 4$;} \\
        (q+1)_2, & \hbox{if $q \equiv 3 \imod 4$.}
  \end{cases}¥
\end{equation}

\begin{thm} \label{thm:constantclassicalbroad}
Let $a$, $q$, $t=t(q)$ be as defined above. Let $n \ge 9$, $G$ a classical group over $\mathbb{F}_{q}$ and let $I(n,t) = [ \frac{1}{2} at \log(n), 4at\log(n))$. Then for a uniformly distributed random $g \in G$, we have
\[\mathbb{P}(|g|_2 \in I(n,t)) \ge \begin{cases}
0.22 & \text{if } G= \GL_n(q), \SL_n(q), \GU_n(q) \text{ or } \SU_n(q);\\
0.11 & \text{if } G= \Sp_{2n}(q), \SO^{\epsilon}_{2n+\delta}(q) \, (\delta \in \{0,1\}).
\end{cases}\]
\end{thm}
\begin{remark} 
 \emph{We note that there are precisely three powers of $2$ in $I(n,t)$. We give positive constant lower bounds  for individual $2$-part orders in Theorem \ref{thm:constantclassical}}.
\end{remark}

  In Theorem \ref{thm:main}, we give lower bounds for a wide range of possible \tpos in classical groups of odd characteristic; however  our methods did not produce effective lower bounds for \emph{all} possible \tpos \hspace{-.11cm}. We give the lower bounds in Theorem \ref{thm:main} in terms of the proportion $p_n(2^{\j})$ of elements in $S_n$ of \tpo $2^{\j}$. Now  Erd{\H{o}}s and Tur\'{a}n \cite{ET2} proved that for a prime power $p_0 \le n$, the proportion $s_{\neg p_0} (n)$ of elements in $S_n$ whose order is not divisible by $p_0$ is
  \begin{align} \label{eqn:erdosturan}
s_{\neg p_0} (n) = \prod_{u=1}^{\lfloor{n/p_0}\rfloor} \left(1- \frac{1}{up_0}\right).
\end{align}
Thus, the proportion $p_n(2^{\j})$ of elements in $S_n$ of \tpo $2^{\j}$ can be evaluated explicitly using the formula
 \begin{align} \label{eqn:pn2j}
p_n(2^{\j}) = s_{\neg 2^{\j+1}}(n) - s_{\neg 2^{\j}}(n).
  \end{align}
 
\begin{thm} \label{thm:main}
Suppose that $\GF$ is a finite classical group defined over a field of odd order $q$, and $t$ is defined by \eqref{eqn:introdeft}.
Let $p_n(2^{\j})$ be defined as above.  Then the proportion of elements in $G$ of \tpo $2^{\a}$ is at least
\begin{small}
\begin{displaymath}
\begin{array}{ll}
 \frac{1}{2\sqrt{2 \pi n}}, & \hbox{if $\GF =\! \GL_{n}(q)$ and $2^{\a}\!=\!(q\!-\!1)_2$, or if $\GF=\GU_n(q)$  and $2^{\a}\!=\!(q\!+\!1)_2$;} \\
\frac{p_n(2^{\a}/t)}{2}, & \hbox{if $\GF = \GL_{n}(q)$ or $\GU_n(q)$, and   $2 \le 2^{\a}/t \le n$;} \\
 \frac{p_{\l}(2^{\a}/t)}{4} , & \hbox{if $\GF = \mathrm{Sp}_{2\l}(q)$,  or $\SO_{2\l+1}(q)$, and $1 \le 2^{\a}/t\le \l$;}\\
\frac{p_{\l}(2^{\a}/t)}{4} , & \hbox{if $\GF = \SO_{2\l}^{\epsilon}(q)$, $1 \le 2^{\a}/t< \l$ and $\l_2 \ne 2^\a/t$;}\\
     \frac{p_{\l}(2^{\a}/t)}{4} - \frac{1}{4\l},           & \hbox{if  $\GF = \SO^{\epsilon}_{2\l}(q)$, $1 \le 2^{\a}/t< \l$, and $\l_2=2^\a/t$.}
 \end{array}
\end{displaymath}
\end{small}
If $2^\a/t \ge 2$ and $2^\a/t$ does not divide $\l$, then the second, third, and fourth lower bounds also hold for $G/Z(G)$. In the first case, a lower bound of $\frac{1}{2\sqrt{2 \pi n}} - \frac{1}{2n}$ holds for $G/Z(G)$.
Furthermore, the lower bounds for $G=\GL_n(q)$ and $\GU_n(q)$ also hold for $G=\SL_n(q)$ and $\SU_n(q)$, unless $2^{\a}/t=n$.
\end{thm}

We remark that Theorems \ref{thm:main} and  \ref{thm:constantSn}(b) imply that for every classical group $G$ in Theorem \ref{thm:main},  there exists an integer $b$ such that the proportion of elements in $G$ with \tpo $2^b$ is at least $ \frac{\sqrt{5}-2}{2}$ when $G = \GL^{\epsilon}_d(q)$ and at least $ \frac{\sqrt{5}-2}{4}$ when $G = \Sp_{d}(q)$ or $\SO^{\epsilon}_d(q)$. Moreover, in Theorem \ref{thm:constantclassical}, we get very close to these lower bounds with a specific power of $2$.

As an application of these results, consider the following problem: Let $C = Cl_m(q) \times Cl_{d-m}(q)$ be the direct product of finite classical groups. For example, $C$ might be the centralizer of an involution contained in a larger classical group $Cl_d(q)$. If we choose an element $(x,y) \in C$ at random, determine the probability that there exists an integer $k$ such that $(x,y)^k = (w,1)$, where $w$ is an involution. Equivalently, find the probability that the \tpo of $x$ is strictly greater than the \tpo of $y$. This problem arises in the analysis of certain black box recognition algorithms, and in particular, the algorithm of Yal\c cinkaya  \cite{Yal1}. Corollary \ref{cor:sukru} shows that this probability is bounded below by $K/d^{1/2}$, where $K$ is an explicit positive constant.
\begin{cor} \label{cor:sukru}
Suppose that $q$ is an odd prime power and that $\GF$ is a finite classical group of dimension $d$ defined over $\mathbb{F}_q$. Consider a direct product of classical groups $C=Cl_m(q) \times Cl_{d-m}(q) \le \GF$, where $2 \le m \le d/2$.  If $(x,y)$ is chosen at random in $C$, then the probability that $|x|_2 > |y|_2$ is at least
\[K/\sqrt{d},\]
where $K$ is an explicit constant, independent of $d$ and $m$.
\end{cor}

The lower bound in Corollary \ref{cor:sukru} applies for all $m \ge 2$; if we specify $m$ greater than $2$, then we obtain much better bounds, which are described in Theorem \ref{thm:sukru}. Furthermore, if we require that $m$ satisfies $d/3\le m \le d/2$, so for example if $C$ is the centralizer of a so-called strong involution (as studied in \cite{LGOB}), then the probability is bounded below by an explicit positive constant. We prove this in Proposition \ref{pro:balanced}. Theorem \ref{thm:oddandtwiceodd} follows from Propositions \ref{caseS2part2} and \ref{pro:caseO2}.  Theorem \ref{thm:constantclassicalbroad} follows from Theorem \ref{thm:constantclassical} and Theorem \ref{thm:main} follows from Propositions \ref{pro:gln}, \ref{pro:caseS2parts2j}, and \ref{pro:caseO2jp}, and Remark \ref{rem:SL}. 
Corollary \ref{cor:sukru} is proved in Section \ref{sec:sukru}.

\section{Preliminaries}

\subsection{Notation}

Let $\Gbar$ be a connected reductive algebraic group defined over $\bar{\mathbb{F}}_q$, the algebraic closure of a field $\mathbb{F}_q$ of odd order $q$.  Let $F$ be a Frobenius morphism of $\Gbar$, and let $G=\Gbar^F$ be the subgroup of $\Gbar$ fixed elementwise by $F$, so that $G$ is a finite group of Lie type. Let $\Tbar$ be an $F$-stable maximal torus in $\Gbar$ and let $W:=N_{\Gbar}(\Tbar)/\Tbar$ denote the Weyl group of $\Gbar$. We will say that two elements $w, w^{\prime}$ in $W$ are $F$-conjugate if there exists $x$ in $W$  such that $w^{\prime}= x^{-1}wF(x)$. This is an equivalence relation, and we will refer to the equivalence classes as \textit{$F$-conjugacy classes}. Moreover, there is an explicit one-to-one correspondence between the $F$-classes of $W$ and the $\GF$-classes of maximal tori in $\GF$. A thorough description of the correspondence can be found in \cite{NiePra}, and we will summarize the necessary results in each section. Although the correspondence is between $F$-classes in $W$ and $\GF$-classes of maximal tori in $\GF$, we will frequently refer to an \textit{element} of $W$ corresponding to a maximal torus in $\GF$, rather than a $\GF$-class of tori, where there is little possibility of confusion.

Now recall that every element $g$ in $\GF$ can be expressed uniquely in the form $g=su$, where $s \in \GF$ is semisimple, $u \in \GF$ is unipotent and $su=us$. This is the multiplicative Jordan decomposition of $g$ (see \cite[p. 11]{Carter2}). We now define a quokka set as a subset of a finite group of Lie type that satisfies certain closure properties.

\begin{dfn}
Suppose that $\GF$ is a finite group of Lie type. A nonempty subset $Q$ of $\GF$ is called a \textit{quokka set}, or quokka subset of $\GF$,  if the following two conditions hold.
\begin{enumerate}
 \item[(i)] For each $g \in \GF$ with Jordan decomposition $g=su=us$, where $s$ is the semisimple part of $g$ and $u$ the unipotent part of $g$, the element $g$ is contained in $Q$ if and only if $s$ is contained in $Q$; and
\item [(ii)] the set $Q$ is a union of $\GF$-conjugacy classes.
\end{enumerate}
\end{dfn}
More generally, one can define a quokka set of an arbitrary finite group (see \cite{NiePra}), but in this paper we will deal exclusively with subsets of finite groups of Lie type. In particular,  we define the subset
\begin{displaymath}
Q(2^{\a}) := \{ g \in \GF \,: \, |g|_{2}=2^{\a} \},
\end{displaymath}
consisting of all the elements $g$ in $\GF$ of \tpo $2^{\a}$.  We readily see that $Q(2^{\a})$ is a quokka set.



Theorem \ref{thm:maintool} will be our main tool to estimate $\frac{|Q(2^{\a})|}{|\GF|}$. It is a direct application of  a theorem of Niemeyer and the second author \cite[Theorem 1.3]{NiePra}.

\begin{thm} \label{thm:maintool}
Let $\GF$ be a finite group of Lie type with Weyl group $W$. Let $Q(2^{\a}):= \{ g \in \GF \,:\, |g|_2=2^{\a} \}$ be the quokka subset of $G$ consisting of the elements of $\GF$ of \tpo $2^{\a}$. Let $C$ be a subset of $W$ consisting of a union of $F$-classes of $W$. For each $F$-class $B$ in $W$, let $T_B$ denote a maximal torus corresponding to $B$.
Then
\begin{align} \label{eqn:maintool1}
\frac{|Q(2^{\a})|}{|\GF|} &\ge \sum_{B \subset C}\frac{|B|}{|W|}.\frac{|T_B \cap Q(2^{\a})|}{|T_B|},
\end{align}
and if there exists a constant $A \in [0,1]$ such that $\frac{|T_B \cap Q(2^{\a})|}{|T_B|} \ge A$ for all $F$-classes in $C$, then
\begin{align} \label{eqn:maintool2}
\frac{|Q(2^{\a})|}{|\GF|} \ge\frac{A|C|}{|W|}.
\end{align}
\end{thm}
\begin{proof}
See \cite[Theorem 1.3]{NiePra}.
\end{proof}
\vspace*{2pt}


To obtain results for some of the projective groups we appeal to Theorem \ref{thm:projectivequokkasets}, the proof of which can be found in \cite[Theorem 1.5]{NiePra}.

\begin{thm}\label{thm:projectivequokkasets}
Let $\GF$ be a finite group of Lie type defined over a field of odd characteristic, let $Z$ be the centre of $\GF$, and let $Q$ be a quokka subset of $\GF$. Then
\[ \frac{|Q|}{|\GF|} \le \frac{|QZ/Z|}{|\GF/Z|} \le \frac{|Q| \cdot |Z|}{|\GF|}.\]
\end{thm}


In light of Lemma \ref{lem:QmeetT} below, we will find lower bounds on $|Q(2^{\a})|/|\GF|$ using Theorem \ref{thm:maintool}, where the subset $C \subset W$ consists of $F$-classes that correspond to classes of maximal tori $T$ whose exponent has $2$-part $2^{\a}$. We refer to such maximal tori $T$ as having \textit{\tpe} \hspace{-0.28cm} $2^{\j}$, and we will write $\exp(T)_2=2^{\a}$ to denote this.

\begin{lemma} \label{lem:QmeetT}
Let $G$ be a finite group and for each $\a$, let $Q(2^{\a})$ be the subset of  $G$ consisting of elements of \tpo $2^{\a}$. Let $T$ be a subgroup of $G$.
\begin{enumerate}
  \item[(a)] If $T$ is cyclic and the $2$-part of $|T|$ is $|T|_2=2^b$, then
\begin{equation} \label{q11}
\frac{|Q(1) \cap T|}{|T|} = \frac{1}{2^{b}};
\end{equation}
  if $ 1  \le j \le b$, then
\begin{equation} \label{q12a}
\frac{|Q(2^{j}) \cap T|}{|T|} = \frac{1}{2^{b+1-j}}.
\end{equation}
  \item[(b)] If $T  = \prod_{i=1}^k C_i$ is a direct product of $k$ cyclic groups $C_i$, all of even order, and $\exp(T)_2=2$, then
\begin{align} \label{q11b1}
\frac{|Q(1) \cap T|}{|T|} = \frac{1}{2^{k}} \quad &\hbox{and} \quad
\frac{|Q(2) \cap T|}{|T|} = \frac{2^k-1}{2^k}.
\end{align}
  \item[(c)] If $T  = \prod_{i=1}^k C_i$ is a direct product of $k$ cyclic groups $C_i$, and $\exp(T)_2=2^b$, then
\begin{equation} \label{q12ak}
\frac{|Q(2^b) \cap T|}{|T|} \ge \frac{1}{2}.
\end{equation}
\item[(d)] Suppose that $T  = \prod_{i=1}^k C_i$ is a direct product of $k$ cyclic groups $C_i$, and that for all $i=1,\ldots, k$, $\exp(C_i)_2=2^b$. Let $z_i$ be the unique involution in $C_i$. Then $1/2^{k}$ of the elements in $T$ have \tpo $2^b$ and satisfy $t^{\frac{|t|}{2}} = (z_1,\ldots,z_k)$.
\end{enumerate}
\end{lemma}

\begin{proof}
(a) Let $x$ be a generator of $T$ so that $T=\langle x \rangle$. Now $|T| = m=2^b c$, for some odd integer $c$, and
\begin{displaymath}
T= \{x^i \,:\, i=1,2, \ldots, 2^bc \}.
\end{displaymath}
The order of an element $x^i$ is $\frac{m}{\gcd(i,m)}$.
We will use induction on the order of $T$. For the base case, we can take the cyclic group of order $2$, and the equations are trivially true for this group. Now suppose that $T$ has order $m=2^bc$ and  $b=0$; so all of the elements of $T$ have odd order and (a) holds trivially, so we may assume that $b \ge 1$.
The elements for which the order has $2$-part equal to $2^b$ are the elements of the form $x^{i}$, where $i$ is odd. The number of odd integers between $1$ and $m=2^bd$ is precisely $m/2$, which proves equation \eqref{q12a} when $j=b$.  The elements of the form $x^i$ where $i$ is even are contained in $\langle x^2 \rangle$, and by induction
\begin{displaymath}
|Q(1) \cap T | = |Q(1) \cap \langle x^2\rangle| = |\langle x^2 \rangle| \frac{1}{2^{b-1}} = |T|\frac{1}{2^b}
\end{displaymath}
thus equation \eqref{q11} holds;  if $b\ge 1$ and $1 \le j \le b-1$, then
\begin{displaymath}
|Q(2^{j}) \cap T|  = |Q(2^{j}) \cap \langle x^2\rangle| = |\langle x^2 \rangle| \frac{1}{2^{b-j}} = |T| \frac{1}{2^{b+1-j}},
\end{displaymath}
proving equation \eqref{q12a} as well and the proof by induction is complete. \\
\noindent
(b) If $\exp(T)_2=2$ and each of the cyclic groups $C_i$ has even order, then an element $t= (c_1, c_2, \ldots, c_k)$ in $T$, where $c_i \in C_i$, has odd order if and only if every $c_i$ has odd order. Thus the first equation in \eqref{q11b1} holds, and the rest of the elements in $T$ have \tpos equal to $2$. So the second equation in \eqref{q11b1} follows as well.\\
\noindent
(c) If $\exp(T)_2=2^b$, then $|C_i|_2 \le 2^b$ for all $i$, and without loss of generality, we may assume that $|C_1|_2=2^b$. If $T$ is cyclic then we appeal to (a), so let us assume that there are at least $2$ direct factors; that is, $k \ge 2$. For this case, we have $T=C_1 \times D$, where $D=\prod_{i=2}^k C_i$ is a product of $k-1$ cyclic groups, and $C_1$ has \tpe $2^b$. From (a), half of the elements of $C_1$ have 2-part order $2^b$. Now observe that if $c \in C_1$ has \tpo $2^b$, then $(c,d) \in C \times D$ also has \tpo $2^b$ for all $d \in D$.  Thus at least half of the elements in $T$ have \tpo $2^b$ also when $k \ge 2$.

\noindent
(d) An element $t \in T$ satisfies the two conditions if and only if $t=(c_1, \ldots, c_k)$ and for all $i=1,\ldots,k$, $c_i$ has \tpo $2^b$. The result follows from \eqref{q12a}.
\end{proof}

 A maximal torus $T$, in a classical group that we will be concerned with, is a direct product of cyclic groups of orders $q^i \pm 1$, for various $i$; so we can use Lemma \ref{lem:2parts} to evaluate $\exp(T)_2$.
\begin{lemma} \label{lem:2parts}
Let $q$ and $i$ be positive integers, and suppose that $q$ is odd. Then
 \begin{equation} \label{qi+1}
  (q^i+1)_2 =
  \left\{
               \begin{array}{ll}
                 2 & \hbox{if $i$ is even;} \\
                 (q+1)_2 & \hbox{if $i$ is odd;}
               \end{array}
             \right.
   \end{equation}
  and
 \begin{equation} \label{qi-1}
  (q^i-1)_2 =
  \left\{
               \begin{array}{ll}
                 (q-1)_2    & \hbox{if $i$ is odd;} \\
                 i_2(q-1)_2 & \hbox{if $q \equiv 1 \imod 4$, and $i$ is even;} \\
                 i_2(q+1)_2 & \hbox{if $q \equiv 3 \imod 4$, and $i$ is even.}
               \end{array}
             \right.
\end{equation}
\end{lemma}
\begin{proof}
First note that  $(q^i-1) = (q-1)(q^{i-1}+q^{i-2} +\cdots+ q+1)$ and if $i$ is odd then $(q^{i-1}+q^{i-2} +\cdots+ q+1)$ is odd and the first line of \eqref{qi-1} follows. The rest of the Lemma is proved in \cite[Lemma 4.1]{LNP}.
\end{proof}

In light of Theorem \ref{thm:maintool}, the proportion of elements of certain \tpos in classical groups is related to the proportion of elements of certain \tpos in $S_n$. We state the main results of \cite{GPsymmetric} on such proportions in $S_n$ that are relevant to our work in this paper. We note that part (c) for $b=0$ was proved in \cite{BLNPS}. 

\begin{thm} \label{thm:constantSn}  \label{t:bestposSn}
Let $a$ and $p_n(2^{\j})$ be as in \eqref{e:a} and \eqref{eqn:pn2j}. \\
\emph{(a)} If $1 \le 2^{\j}=\alpha\log(n) \le n$ and $\alpha$ is contained in one of the intervals $I$ in column 1 of Table  \ref{tab:forpropenot2e}, and $n \ge N_I$, where $N_I$ is an integer given in column 2 of Table  \ref{tab:forpropenot2e}, then a constant lower bound for $p_n(\alpha \log(n))$ is given in column 3 of Table \ref{tab:forpropenot2e}. \\
\emph{(b)} For all positive integers $n$,  there exists  a positive integer $b$  such that $p_n(2^b) \ge \sqrt{5}-2$.\\
\emph{(c)}  If $1 \le 2^{\j}\le n$, then there exists  a constant $K_{\j}$ depending only on $\j$ such that $p_n(2^{\j}) \ge K_{\j} /n^{1/2^{\j+1}}$ and $K_0 = (2 \pi)^{-1/2}$.
\end{thm}
%

%

\section{A more detailed version of Theorem \ref{thm:constantclassicalbroad}}
Note that  $a$ and $t=t(q)$ are as defined in \eqref{e:a} and \eqref{eqn:introdeft} and that $q$ is an \emph{odd} prime power.
Theorem \ref{thm:constantclassical} below follows from Theorems \ref{thm:main} and  \ref{thm:constantSn}, and Remarks \ref{lin:semiregular}, \ref{rem:6.4}, and \ref{rem:7.3}. 

\begin{thm} \label{thm:constantclassical}
 Suppose that $n \ge N'_I$, where $N'_I$ is an integer given in column 4 of Table \ref{tab:forpropenot2e}, $\alpha \log(n)$ is a power of $2$, with $\alpha$ contained in one of the intervals $I$ in column 1 of Table \ref{tab:forpropenot2e}, and let $2^\j= \alpha \log(n)t$ (so in all cases $t \le 2^{\j} \le nt$). Let $\GF = \GL_n(q)$, $\SL_n(q)$, $\GU_n(q)$, $\SU_n(q)$, $\Sp_{2n}(q)$, $\SO_{2n+1}(q)$, or $\SO^{\epsilon}_{2n}(q)$. Then columns 5 and 6 of Table \ref{tab:forpropenot2e} give explicit constant lower bounds on the proportion of elements in $\GF$, and in $\GF/Z(\GF)$, with \tpo $2^{\j}$. Moreover, when $N'_I> n\ge N_I$, the lower bounds in column 6 also hold for this proportion in $G$ (but not necessarily in $G/Z(G)$).
\end{thm}


{\openup 8pt
\fontsize{8}{16}
\selectfont
\begin{table}
\setlength{\tabcolsep}{3mm}
   \caption{Lower bounds on the proportion of elements in $G$ or $G/Z(G)$ whose \tpo is $2^{\j}$.}\label{tab:forpropenot2e}
\vspace{1mm}
    \begin{center}
    \begin{tabular}{cccccc}
\hline
\multicolumn{1}{c}{Interval $I$}  & \multicolumn{1}{c}{$N_I$}&\multicolumn{1}{c}{$G=S_n$, }&\multicolumn{1}{c}{$N'_I$} &\multicolumn{1}{c}{$G=\GL_n(q)$, $\SL_n(q)$,} & \multicolumn{1}{c}{$\GF= \Sp_{2n}(q)$, $\SO^{\pm}_{2 n}(q)$,} \rule{0cm}{3.5mm} \\
\multicolumn{1}{c}{containing $\alpha$}&&\multicolumn{1}{c}{$2^{\j}=\alpha \log(n)$,}&&\multicolumn{1}{c}{$\GU_n(q)$, or $\SU_n(q)$;}& \multicolumn{1}{c}{or $\SO_{2 n+1}(q)$; } \\
&&\multicolumn{1}{c}{$n \ge N_I$}&&\multicolumn{1}{c}{$2^{\j}=\alpha \log(n)t$,}& \multicolumn{1}{c}{$2^{\j}=\alpha \log(n)t$,}\\
&&&&\multicolumn{1}{c}{$n \ge N'_I$}& \multicolumn{1}{c}{ $n \ge N'_I$} \rule[-4pt]{0mm}{3mm} \\
\hline
\rule{0cm}{3mm}$\left[\frac{a}{8}, \frac{a}{4}\right)$ &47& 0.0165662676 &2208
 &0.0082331338
 & 0.0041415669\\
 $\left[\frac{a}{4}, \frac{a}{2}\right)$ \rule{0cm}{3mm} &7& 0.1170040878 &47 &0.0584520439
 &0.0292260219
 \\
 $\left[\frac{a}{2}, a\right)$  \rule{0cm}{3mm}&3& 0.2351203790 &7&  0.1175101895
&0.0587550947
  \\
$[a,2a)$ &2& 0.1531015975 &3 &0.0762123059
 &0.0381061529 \\
$[2a,4a)$ &8&
  0.0605468750 &9 & 0.0292468750 &	0.0146234375
 \\
$[4a,8a)$   &16& 0.0307617187&17  & 0.0090952640	& 0.0045476320
 \\
$[8a,16a)$ &64&   0.0155029296&65  & 0.0036895334 &	0.0018447667   \\ $[16a,32a)$  &128& 0.0065085752 &129  &0.0018617415 &	0.0009308707
\\ \hline
\end{tabular}
 \end{center}
  \end{table}}
Theorem \ref{thm:constantclassicalbroad} follows easily from Theorem \ref{thm:constantclassical}, which we now demonstrate.
\begin{proof}[Proof of Theorem \ref{thm:constantclassicalbroad}]
  This follows easily by adding up the entries in rows $3$, $4$ and $5$ of columns $5$ and $6$ of Table \ref{tab:forpropenot2e}.
\end{proof}

\section{Stirling numbers}
To obtain estimates for the proportion of maximal tori with certain \tpes \hspace{-0.11cm}, it will be useful to know how many permutations there are in $S_{\l}$ with exactly $k$ cycles. These numbers are known as the \emph{unsigned Stirling numbers of the first kind}, and we will denote them by  $c(\l,k)$. The \emph{signed Stirling numbers of the first kind} $s(\l,k)$ are then defined by the equation
\begin{align}
s(\l,k) = (-1)^{\l-k}c(\l,k).
\end{align}
In the sequel, we will need to calculate various summations involving $c(\l,k)$ and $s(\l,k)$; we devote this section to these calculations.

For $n\ge 1$, the signed Stirling numbers of the first kind have generating function
\begin{align} \label{eqn:slkgeneratingfn}
\sum_{k=1}^{\l} s(\l,k)x^k&= \prod_{i=0}^{\l-1} (x-i),
\end{align}
which is valid for $-1<x<1$ (see \cite[24.1.3]{abramowitz}). Substituting $-x$ into \eqref{eqn:slkgeneratingfn} gives
 \begin{align}
\sum_{k=1}^{\l} s(\l,k)(-1)^k x^k&= (-1)^{\l}\prod_{i=0}^{\l-1} (x+i). \notag
\end{align}
Thus, we obtain the generating function for the unsigned Stirling numbers of the first kind:
\begin{align} \label{eqn:clkgeneratingfn}
\sum_{k=1}^{\l} (-1)^{k-\l}s(\l,k) x^k &=\sum_{k=1}^{\l} c(\l,k) x^k= \prod_{i=0}^{\l-1} (x+i),
\end{align}
which is also defined for $-1<x<1$.  We substitute $x=1/2$ and $x=1/4$ into equations \eqref{eqn:slkgeneratingfn} and \eqref{eqn:clkgeneratingfn}, giving us four equations, and we define these quantities as $A_{\l}$, $B_{\l}$, $C_{\l}$, and $D_{\l}$ respectively, so that we can refer to them later in this section. We have
\begin{small}
\begin{align}
A_{\l} :=\sum_{k=1}^{\l}   \frac{c(\l,k)}{2^{k}} &= \frac{1}{2} \frac{3}{2} \cdots \frac{2 \l -1}{2}  = \frac{(2 \l)!}{2^{2\l} \l !}, \label{eqn:clk2^k}\\
B_{\l} :=\sum_{k=1}^{\l} \frac{s(\l,k)}{2^k}   &=
 \frac{1}{2} \frac{(-1)}{2} \frac{(-3)}{2} \cdots \frac{(3-2 \l )}{2}  = \frac{(-1)^{\l-1}(2 \l-2)!}{2^{2\l-1} (\l-1) !}, \label{eqn:slk2^k}\\
C_{\l} :=\sum_{k=1}^{\l}  \frac{c(\l,k)}{4^{k}} &= 4^{-\l}\prod_{i=1}^{\l} (4i-3)= \frac{\Gamma(\l+1/4)}{\Gamma(1/4)} \hbox { (see \cite[6.1.10]{abramowitz}), and}\label{eqn:clk4^k}\\
D_{\l} := \sum_{k=1}^{\l}\! \frac{s(\l,k)}{4^k}   \!&= \!\frac{(-1)^{\l-1}}{4^{\l}}\prod_{i=1}^{\l\!-\!1} (4i\!-\!1)=\frac{(-1)^{\l \!-\!1}\Gamma(\l\!-\! \frac{1}{4})}{4\Gamma(3/4)} \hbox{ (see \cite[6.1.14]{abramowitz})\!}. \label{eqn:slk4^k}
\end{align}
\end{small}
On the one hand, we have
\begin{align*}
\frac{A_{\l} +(-1)^{\l}B_{\l}}{2} &=\sum_{k=1}^{\l} \frac{c(\l,k)(1+(-1)^{2\l-k})}{2^{k+1}} =\sum_{\substack{{k=1}\\{k \text{ even}}}}^{\l} \frac{c(\l,k)}{2^k},
 \end{align*}
 and on the other hand, we have
\begin{align*}
\frac{A_{\l} +(-1)^{\l}B_{\l}}{2}  &=\frac{1}{2} \left(\frac{(2 \l)!}{2^{2\l} \l !} + \frac{(-1)^{2\l-1}(2 \l-2)!}{2^{2\l-1} (\l-1) !} \right)\\
&= \frac{(2\l)!}{2^{2\l+1}\l!}\left(\frac{2\l-2}{2\l-1}\right);
\end{align*}
thus
\begin{align}
\sum_{\substack{{k=1}\\{k \text{ even}}}}^{\l} \frac{c(\l,k)}{2^k}  &=  \frac{(2\l)!}{2^{2\l+1}\l!}\left(\frac{2\l-2}{2\l-1}\right) \label{eqn:stirling2^keven}.
\end{align}
Similarly,  we have
\begin{align*}
\frac{A_{\l} + (-1)^{\l+1}B_{\l}}{2}  &= \sum_{k=1}^{\l} \frac{c(\l,k)(1+(-1)^{2\l+1-k})}{2^{k+1}} = \sum_{\substack{{k=1}\\{k \text{ odd}}}}^{\l} \frac{c(\l,k)}{2^k},
\end{align*}
 and so
\begin{equation} \label{eqn:stirling2^kodd}
   \begin{split}
 \sum_{\substack{{k=1}\\{k \text{ odd}}}}^{\l} \frac{c(\l,k)}{2^k}= \frac{ A_{\l} + (-1)^{\l+1}B_{\l}}{2}&= \frac{1}{2} \left(\frac{(2 \l)!}{2^{2\l} \l !} + \frac{(-1)^{2\l}(2 \l-2)!}{2^{2\l-1} (\l-1) !} \right)\\
&= \frac{(2\l)!}{2^{2\l+1}\l!}\left(\frac{2\l}{2\l-1}\right). 
\end{split}
\end{equation}
Doing the same for $C_\l$ and $D_\l$ gives the following equations:
\begin{align}
\sum_{\substack{{k=1}\\{k \text{ even}}}}^{\l} \frac{c(\l,k)}{4^k}&= \frac{1}{2}\left( \frac{\Gamma(\l+1/4)}{\Gamma(1/4)} - \frac{\Gamma(\l-1/4)}{4\Gamma(3/4)} \right);
 \label{eqn:stirling4^keven}\\
 \sum_{\substack{{k=1}\\{k \text{ odd}}}}^{\l} \frac{c(\l,k)}{4^k}  &=\frac{1}{2}\left( \frac{\Gamma(\l+1/4)}{\Gamma(1/4)} + \frac{\Gamma(\l-1/4)}{4\Gamma(3/4)} \right) \label{eqn:stirling4^kodd}.
\end{align}
\noindent If $\l \ge 2$, then Stirling's approximation (see \cite{Stirlingapprox}) implies
\begin{displaymath}
\left| \frac{\l!}{\sqrt{2\pi \l}\l^\l e^{-\l}} -1- \frac{1}{12\l} \right| \le \frac{1}{288\l^2} + \frac{1}{9940\l^3},
\end{displaymath}
and this enables us to show that if $\l \ge 2$, then
 \begin{align} \label{eqn:2ellfactorial}
\frac{(2 \l)!}{(\l !)^2 2^{2\l}} &\ge \frac{25}{29 \sqrt{\pi\l}}.
\end{align}
Moreover, in \cite{vasicgammabds}, Ke{\v{c}}ki{\'c} and Vasi{\'c} prove that if $y> x \ge 1$, then
 \begin{align} \label{eqn:gammaquotientbounds}
\frac{x^{x-1/2}}{y^{y-1/2}}e^{y-x}< \frac{\Gamma(x)}{\Gamma(y)}&< \frac{x^{x-1}}{y^{y-1}}e^{y-x}.
 \end{align}
 Recall that $\Gamma(\l+1) = \l !$. We can use \eqref{eqn:gammaquotientbounds} to show that
\begin{align*}
 e^{3/4} \frac{\left(1-\frac{3}{4(\l+1)}\right)^{\l-1/4}}{\Gamma(1/4)(\l+1)^{3/4}} <\sum_{k=1}^{\l} \frac{c(\l,k)}{4^{k}\l!} &< \frac{e^{3/4}}{(\l+1)^{3/4}\Gamma(1/4)}.
\end{align*}
Now  $\left(1-\frac{3}{4(\l+1)}\right)^{\l-1/4} > \left(1-\frac{3}{4(\l+1)}\right)^{\l}$ and it is elementary to show that  $\left(1-\frac{3}{4(\l+1)}\right)^{\l}>e^{-3/4}$ (using \cite[Lemma 2.3]{GPsymmetric} for example); thus we have
\begin{align*}
\frac{1}{\Gamma(1/4)(\l+1)^{3/4}}<\sum_{k=1}^{\l} \frac{c(\l,k)}{4^{k}\l!} =\frac{\Gamma(\l+1/4)}{\Gamma(\l+1)\Gamma(1/4)}&< \frac{e^{3/4}}{(\l+1)^{3/4}\Gamma(1/4)}.
\end{align*}
It follows that
\begin{align} \label{eqn:clk4^kasl^3/4}
\frac{1}{4(\l+1)^{3/4}}<\sum_{k=1}^{\l} \frac{c(\l,k)}{4^{k}\l!}=\frac{\Gamma(\l+1/4)}{\Gamma(\l+1)\Gamma(1/4)} < \frac{3}{5(\l+1)^{3/4}}.
\end{align}
Similarly, to approximate the quantities in \eqref{eqn:stirling4^keven} and \eqref{eqn:stirling4^kodd}, observe that, using \eqref{eqn:gammaquotientbounds},
\begin{align*}
\frac{1}{4\Gamma(3/4)(\l+1)^{1/4}}< \frac{\Gamma(\l-1/4)}{4\Gamma(3/4)\l!} &< \frac{e^{5/4}}{4\Gamma(3/4)(\l+1)^{5/4}},
\end{align*}
and so we have
\begin{align}
\frac{1}{5(\l+1)^{5/4}}< \frac{\Gamma(\l-1/4)}{4\Gamma(3/4)\l!} &< \frac{18}{25(\l+1)^{5/4}}.
\end{align}
Thus, we have the following bounds:
\begin{small}
\begin{align}
\frac{1}{8(\l+1)^{3/4}} \!-\!\frac{9}{25(\l+1)^{5/4}}<\sum_{\substack{{k=1}\\{k \text{ even}}}}^{\l} \frac{c(\l,k)}{\l!4^k}&<\frac{3}{10(\l+1)^{3/4}} \!-\! \frac{1}{10(\l+1)^{5/4}}; \label{eqn:boundsclk4^keven}\\
\frac{1}{8(\l+1)^{3/4}}\!+\!\frac{1}{10(\l+1)^{5/4}}<\sum_{\substack{{k=1}\\{k \text{ odd}}}}^{\l} \frac{c(\l,k)}{\l!4^k}&<\frac{3}{10(\l+1)^{3/4}} \!+\! \frac{9}{25(\l+1)^{5/4}}. \label{eqn:boundsclk4^kodd}
\end{align}
\end{small}
Finally, we will also need to know how many permutations in $S_\l$ have $k$ cycles, an even number of them having even length. We say that such a permutation has $k$ cycles and an even number of even cycles.
\begin{lemma} \label{lem:numberevencycles}
Let $c_1(\l,k)$ be the number of permutations in $S_{\l}$ with $k$ cycles, and an even number of even cycles. Let $c_2(\l,k)$ be the number of permutations in $S_{\l}$ with $k$ cycles, and an odd number of even cycles. Then
\begin{displaymath}
c_1(\l,k) = \left\{
  \begin{array}{ll}
    c(\l,k), & \hbox{if $\l$ and $k$ are either both even, or both odd;} \\
    0, & \hbox{otherwise;}
  \end{array}
\right.
\end{displaymath}
and
\begin{displaymath}
c_2(\l,k) = \left\{
  \begin{array}{ll}
    0, & \hbox{if $\l$ and $k$ are either both even, or both odd;} \\
    c(\l,k), & \hbox{otherwise.}
  \end{array}
\right.
\end{displaymath}
\end{lemma}
\begin{proof}
Consider a conjugacy class in $S_{\l}$ corresponding to the partition $(d_1, \ldots, d_k)$. In particular, these permutations are some of the $c(\l,k)$ permutations in $S_{\l}$ with $k$ cycles, where $c(\l,k)$ is an unsigned Stirling number of the first kind. Suppose that $m$ of the $d_i$ are even. We have
\begin{displaymath}
\l = d_1 + \cdots + d_k
\end{displaymath}
and reducing this equation modulo 2 gives
\begin{displaymath}
\l \equiv k-m \imod 2,
\end{displaymath}
or equivalently,
\begin{displaymath}
m \equiv k-\l \imod 2.
\end{displaymath}
Thus there are an even number $m$ of even cycles if and only if $k$ and $\l$ are both even or both odd. Therefore the first equation holds. The second equation follows from the first since
\[
c_2(\l,k) =  c(\l,k) - c_1(\l,k). \qedhere
\]
\end{proof}

\section{Proof of theorem \ref{thm:main} in linear and unitary cases} \label{section:L}

In this section we prove the following proposition, which proves Theorem~\ref{thm:main} in the case where $\GF= \GL_n(q)$ or $\GU_n(q)$. We use the following notation: $\GL^{\epsilon}_n(q)$ refers to $\GL_n(q)$ when $\epsilon=+$, and $\GU_n(q)$ when $\epsilon=-$; for every nonnegative integer $\a$, we define $Q(2^{\a})$ to be the quokka set
\[Q(2^{\a})= \{g \in \GF \,:\, |g|_2=2^{\a}\};\]
we let $p_n(2^{\j})$ denote the proportion of elements in $S_n$ with \tpo $2^{\j}$;
and we let $t=t(q)$ be as in \eqref{eqn:introdeft}. 
Proposition \ref{pro:gln} proves Theorem \ref{thm:main} in the linear and unitary cases.
\begin{pro} \label{pro:gln}
Suppose that $\GF = \GL^{\epsilon}_n(q)$, and $q$ is odd.
If $2 \le 2^\a/t \le n$, then
\begin{align}\label{eqn:lin1}
\frac{|Q(2^\a)|}{|\GF|} \ge   \frac{p_n(2^\a/t)}{2} 
\end{align}
and
\begin{align}\label{eqn:lin2}
\frac{|Q((q-\epsilon)_2)|}{|\GF|} \ge  \frac{1}{2\sqrt{2 \pi n}}.
\end{align}
If $2^b/t$ does not divide $n$, then \eqref{eqn:lin1} holds for $\GF/Z(\GF)$. If $n$ is odd, then \eqref{eqn:lin2} holds for $\GF/Z(\GF)$, and if $n$ is even, then \eqref{eqn:lin2} holds for $G/Z(G)$ if we replace the lower bound with $\frac{1}{2\sqrt{2 \pi n}} - \frac{1}{2n}$.
\end{pro}

\begin{remark} \label{rem:missingorders}
\emph{
Proposition \ref{pro:gln} does not provide any information on the proportion of odd order elements. Moreover, if $q \equiv 1 \imod 4$ and $\epsilon=+$, then $(q-\epsilon)_2=t$, so Proposition \ref{pro:gln} provides lower bounds for the proportion of elements of \tpo $2^{\a}$ for all $2^{\a} \ge t$, but provides no information when $2^{\a}<t$. Similarly, if $q \equiv 3 \imod 4$ and $\epsilon=-$, then $(q-\epsilon)_2=t$ and we have the same situation: we have no information when $2^{\a}<t$. If $q \equiv 1 \imod 4$ and $\epsilon=-$, or if $q \equiv 3 \imod 4$ and $\epsilon=+$, then $(q-\epsilon)_2<t$ and so Proposition \ref{pro:gln} gives lower bounds on $|Q(2^{\a})|/|\GF|$ when $2^{\a}=2$ and when $2^{\a} \ge 2t$, but provides no information if $4 \le 2^{\a} \le t$.}
\end{remark}

\begin{proof}
In Lemma \ref{f=eodd} below, we find the \tpe of a maximal torus corresponding to an element $w \in W=S_n$.

\begin{lemma} \label{f=eodd}
 Let $\GF=\GL^{\epsilon}_n(q)$ and let $T$ be a maximal torus in $\GF$ corresponding to an element $w$ in $W=S_n$. If $2^{\j}$ is the \tpo of $w$, then
\begin{displaymath}
\exp(T)_2 = \left\{
              \begin{array}{ll}
                 2^{\j}t, & \hbox{if  $\j \ge 1$;} \\
                (q-\epsilon)_2, & \hbox{if $\j=0$.}
              \end{array}
            \right.
\end{displaymath}

\end{lemma}
\begin{proof}
 A maximal torus $T$ in $\GL_n(q)$, corresponding to the partition $(a_1, \ldots, a_k)$, is of the form
\begin{displaymath}
T= \prod_{i=1}^k (q^{a_i}-1),
\end{displaymath}
and a maximal torus $T$ in $\GU_n(q)$ is of the form
\begin{displaymath}
T= \prod_{i=1}^k (q^{a_i}-(-1)^{a_i}).
\end{displaymath}
Since $w$ has \tpo $2^\a$, we know that the maximum $2$-part of the $a_i$ is $2^\a$, and Lemma \ref{lem:2parts} implies the result.
\end{proof}


We apply Theorem \ref{thm:maintool}, with $C \subset W$ consisting of classes in $W=S_n$ corresponding to maximal tori of \tpe $2^{\a}$. For any such maximal torus,  we have
\begin{displaymath}
\frac{|Q(2^{\a})\cap T|}{|T|} \ge \frac{1}{2},
\end{displaymath}
 by Lemma \ref{lem:QmeetT}, and so we can take $A=1/2$ in Theorem \ref{thm:maintool}.  Thus, we have
\begin{displaymath}
\frac{|Q(2^{\a})|}{|\GF|} \ge \frac{|C|}{2n!},
\end{displaymath}
where $C$ consists of the classes of permutations in $W=S_n$ with \tpo $2^{\a}/t$; now \eqref{eqn:lin1} follows.

Similarly, if  $\exp(T)_2=(q-\epsilon)_2$, then by Lemma \ref{lem:QmeetT} (b) and (c),
\begin{displaymath}
\frac{|Q((q- \epsilon)_2)\cap T|}{|T|} \ge \frac{1}{2}
\end{displaymath}
and so we again take $A=1/2$ in Theorem \ref{thm:maintool}. Now $T$ has \tpe $(q-\epsilon)_2$ if and only if $T$ corresponds to $w \in S_n$ of odd order (by Lemma \ref{f=eodd}). Thus $C \subset S_n$ just consists of odd order elements, and Theorem \ref{thm:maintool} implies that
\begin{displaymath}
\frac{|Q((q- \epsilon)_2)|}{|\GF|} \ge \frac{p_n(1)}{2},
\end{displaymath}
which is just \eqref{eqn:lin2}, by Theorem \ref{thm:constantSn}(c).

Now suppose that $g \in Q(2^b)$, with $2^{\j} \ge 2$, and $g^{|g|/2} \in Z(G)$. If $g$ is contained in a maximal torus $T$, then $g$ must have \tpo $2^b$ in each direct factor of $T$. In particular, by Lemma \ref{f=eodd}, $T$ must correspond to a partition $(a_1, \ldots ,a_k)$ such that $2^b/t$ divides $a_i$ for every $i$. Thus, by Theorem \ref{thm:projectivequokkasets}, if $2^b/t$ does not divide $n$, then $Q(2^b)Z(G)/Z(G)$ is a quokka set and consists of elements of order $2^b$ in $\GF/Z(\GF)$, and \eqref{eqn:lin1} holds in $\GF/Z(\GF)$.

We now prove \eqref{eqn:lin2} for $G/Z(G)$. Let $Q'((q- \epsilon)_2)$ denote the set of elements $g$ in $G$ with \tpo $(q- \epsilon)_2$ such that $g^{\frac{|g|}{2}} \not\in Z(G)$.   Note that any  maximal torus $T$ corresponding to an odd order permutation must have \tpe $(q- \epsilon)_2$. In particular, each of the $k$ direct factors of $T$ has \tpe $(q- \epsilon)_2$ and  by Lemma \ref{lem:QmeetT}(d), we have
 \begin{displaymath}
\frac{|Q'((q- \epsilon)_2) \cap T|}{|T|} = 1- \frac{1}{2^{k}} - \frac{1}{2^k} = 1- \frac{1}{2^{k-1}}.
 \end{displaymath}
Since this is at least $1/2$ if $k \ge 2$, we take $C \subset S_n$ to consist of odd permutations with at least $2$ cycles. If $n$ is even, then $\frac{|C|}{|W|} = s_{\neg 2}(n)$ as before. If $n$ is odd, then $\frac{|C|}{|W|} = s_{\neg 2}(n)- \frac{1}{n}$, since the proportion of $n$-cycles in $S_n$ is $\frac{1}{n}$. Applying Theorem \ref{thm:maintool} implies that
\begin{align}
\frac{|Q'((q- \epsilon)_2)|}{|G|}& \ge \frac{(s_{\neg 2}(n)- \frac{((2,n)-1)}{n})}{2}.\notag
\end{align}
Now every element in $Q'((q-\epsilon)_2)$ has \tpo $(q- \epsilon)_2$ in $G/Z(G)$, so Theorem \ref{thm:projectivequokkasets} implies the results for $G/Z(G)$.
\end{proof}

\begin{remark} \label{rem:SL}
\emph{When $\GF=\SL^{\epsilon}_n(q)$ and $2^{\a}/t<n$, a maximal torus $T$ has \tpe $2^{\a}$ if and only if the corresponding maximal torus in $\GL^{\epsilon}_n(q)$ has \tpe $2^{\a}$. Thus the proofs for $\GF=\SL^{\epsilon}_n(q)$ apply unless $n=2^{\a}/t$. In this exceptional case, there are no tori in $\SL_n^\epsilon(q)$ with \tpe $2^{\a}$, and no elements in $\GF$ of \tpo $2^{\a}$ in $\SL_n^\epsilon(q)$.}
\end{remark}

\begin{remark} \label{lin:semiregular}
\emph{If $2^b/t \ge 2$ and $2^b/t$ divides $n$, then the proportion of elements in $S_n$ with partitions of the form $(a_1, \ldots, a_k)$ such that every $a_i$ has $2$-part $2^b/t$ is at most
\begin{displaymath}
3\l^{ \frac{t}{2^{b+1}}-1},
\end{displaymath}
by \cite[Theorem 1.1]{NPPY}.  Thus, by only taking the $F$-classes in $C$ that are not of this form, we obtain
\begin{displaymath}
\frac{|Q(2^\a)|}{|\PGL^{\epsilon}_n(q)|} \ge   \frac{p_n(2^\a/t)}{2} - \frac{3}{2\l^{1-\frac{t}{2^{b+1}}}}.
\end{displaymath}
To verify the constant lower bounds for $G/Z(G)$ claimed in Theorem \ref{thm:constantclassical}, we find constant lower bounds on the proportion of elements in $S_n$ whose order is $\alpha\log(\l)$ and does not have cycle structure of the form described above. We can verify this for $\l \le 85$ by summing up the relevant conjugacy class sizes in MAGMA. For $85 \le \l \le 10^6$, we calculate $p_{\l}(2^\a) - 3\l^{\frac{t}{2^{b+1}}-1}$ explicitly in MAGMA using \eqref{eqn:pn2j}. For $\l > 10^6$, we use the lower bounds on $p_{\l}(2^\a)$ from Theorem \ref{thm:constantSn} together with the upper bound
\[3\l^{ \frac{t}{2^{b+1}}-1} \le \frac{3}{10^{6(3/4)}}= \frac{3}{10^{9/2}}. \]
This remark completes the proof of Theorem \ref{thm:constantclassical} in the linear and unitary cases.}
\end{remark}

\section{Symplectic and odd dimensional orthogonal cases}

Suppose that $\GF = \mathrm{Sp}_{2\l}(q)$. Then $W=C_2 \wr S_{\l} \le S_{2 \l}$, and a partition
\[\beta = (\beta^+, \beta^-)=(b_1^+,\ldots, b_m^+, b_1^-, \ldots, b_l^-)\]
 of $\l$, representing an $F$-class in $W$, corresponds to a $\GF$-class of maximal tori $T$ of the form
\begin{equation}
T= \prod_{i=1}^m (q^{b_i^+}-1) \times \prod_{j=1}^l (q^{b_j^-}+1).
\end{equation}
For $w \in W=C_2 \wr S_{\l}$, write $\sigma(w)$ for the natural image of $w$ in $S_{\l}$. The partition $\beta$, without the signs, represents the cycle structure of $\sigma(w)$ in $S_{\l}$; the $+$ sign indicates that the $b_1^{+}$ cycle in $S_{\l}$ is the image of two cycles of length $b_1$ in $W\le S_{2\l}$, and these will be referred to as \textit{positive cycles}. The $-$ sign indicates that the $b_1^{-}$ cycle in $S_{\l}$ is the image of one cycle of length $2b_1$ in $W\le S_{2\l}$, and  these will be referred to as \textit{negative cycles}.

If $\GF = \SO_{2\l+1}(q)$, then the Weyl group $W$ and the correspondence between classes in $W$ and $\GF$-classes of maximal tori is exactly the same. So without loss of generality, we may assume that $\GF= \mathrm{Sp}_{2\l}(q)$.

Using Lemma \ref{lem:2parts}, we can calculate the \tpe for each maximal torus in $\GF$.
\begin{lemma}\label{caseS2parts}
Let $\GF = \mathrm{Sp}_{2\l}(q)$ or $\SO_{2\l+1}(q)$ and let $T$ be a maximal torus in $\GF$ corresponding to the partition $\beta = (\beta^+, \beta^-)=(b_1^+,\ldots, b_m^+, b_1^-, \ldots, b_l^-)$ of $\l$. If $q \equiv 1 \imod 4$, then
\begin{displaymath}
\exp(T)_2=\left\{
  \begin{array}{ll}
    (q-1)_2 \max_i \{(b_i^+)_2\}, & \hbox{if there is at least one positive cycle;} \\
    2, & \hbox{if there are no positive cycles;}
  \end{array}
\right.
\end{displaymath}
and if $q \equiv 3 \imod 4$, then
\begin{displaymath}
\exp(T)_2=\left\{
  \begin{array}{ll}
    (q+1)_2 \max_i \{(b_i^+)_2\}, & \hbox{if at least one of the $b_i^+$ is even;} \\
    (q+1)_2, & \hbox{if none of the $b_i^+$ are even, but at least} \\
             &  \hbox{one of the $b_j^-$ is odd;} \\
    2, & \hbox{otherwise, that is if all of the $b_i^+$ are odd} \\
        & \hbox{and all of the $b_j^-$ are even.}
  \end{array}
\right.
\end{displaymath}
\end{lemma}
\begin{proof}
Since $\exp(T)_2$ is just the maximum $2$-part of all of the $q^{b_i^+}-1$ and $q^{b_j^-}+1$, the result follows from Lemma \ref{lem:2parts}.
\end{proof}

So let us first investigate the proportion of elements with \tpos $1$ and $2$. Proposition \ref{caseS2part2} proves Theorem \ref{thm:oddandtwiceodd} in the cases where $G$ is symplectic or odd dimensional orthogonal.
\begin{pro} \label{caseS2part2}
Let $\l \ge 1$ and let $(\GF,\delta_1,\delta_2) = (\Sp_{2\l}(q),1,1)$, $(\SO_{2\l+1}(q),1,1)$, $(\Om_{2\l+1}(q),2,2)$, $(\PSp_{2\l}(q),1,2)$, $(\PSO_{2\l+1}(q),1,2)$, or $(\POm_{2\l+1}(q),2,4)$.
If $Q(1):= \{g\in \GF \,:\, |g|_2=1\}$, then the proportion of odd order elements in $\GF$ is
\begin{align}
\frac{|Q(1)|}{|\GF|} \ge \frac{\Gamma(\l+1/4)\delta_1}{\Gamma(1/4)\l!}  \ge \frac{1}{4 (\l+1)^{3/4}}. \label{eqn:Sodd}
\end{align}
Also, if $Q(2):= \{g\in \GF \,:\, |g|_2=2\}$, then the proportion of elements of $\GF$ with $2$-part order precisely $2$ is
\begin{align}
\frac{|Q(2)|}{|\GF|} \ge \frac{(2 \l)!}{2^{2\l} (\l !)^2} - \frac{\delta_2\Gamma(\l+1/4)}{\Gamma(1/4)\l!} \ge \frac{25}{29 \sqrt{\pi \l}}-\frac{3\delta_2}{5 (\l+1)^{3/4}}. \label{eqn:S2part2}
\end{align}
\end{pro}

\begin{proof}
To prove both inequalities, we will  apply Theorem \ref{thm:maintool} with $C \subset W$ as the union of $F$-classes in $W$ that correspond to maximal tori of \tpe 2.  To prove \eqref{eqn:S2part2},  first suppose that $q \equiv 1 \imod 4$ so that, by Lemma \ref{caseS2parts},  $C$ consists of $F$-classes whose partitions $\beta=(b_1^+,\ldots, b_m^+, b_1^-, \ldots, b_l^-)$ have no positive cycles; that is, $m=0$. Now every permutation $\tau \in S_{\l}$ with $k$ cycles corresponds to precisely ${2^{\l-k}}$ elements $w$ in $W$ that satisfy the conditions: $\sigma(w)=\tau$, and each of the $k$ cycles is negative. Moreover, given a maximal torus  with \tpe $2$ and $k$ direct factors, Lemma \ref{lem:QmeetT} implies that
\begin{displaymath}
\frac{|Q(2) \cap T |}{|T|} = \frac{2^k-1}{2^{k}}.
\end{displaymath}
 Thus, by Theorem \ref{thm:maintool}, we have
\begin{align}
\frac{|Q(2)|}{|\GF|} &\ge \sum_{k=1}^{\l} \frac{c(\l,k)}{2^k\l!}  \left(\frac{2^k-1}{2^{k}} \right) =   \sum_{k=1}^{\l} \frac{c(\l,k)}{2^k\l!}  - \sum_{k=1}^{\l} \frac{c(\l,k)}{4^k\l!} . \label{eqn:6.4}
\end{align}
Applying \eqref{eqn:clk2^k} and \eqref{eqn:clk4^k} to \eqref{eqn:6.4} gives
\begin{displaymath}
\frac{|Q(2)|}{|\GF|}  \ge \frac{(2 \l)!}{2^{2\l} (\l !)^2} -   \frac{\Gamma(\l+1/4)}{\Gamma(1/4)\l!},
\end{displaymath}
which is the first inequality claimed in \eqref{eqn:S2part2}. The second inequality in \eqref{eqn:S2part2} now follows from \eqref{eqn:2ellfactorial} and \eqref{eqn:clk4^kasl^3/4} for $\l \ge 2$, and an explicit calculation verifies the case that $\l=1$.
This proves \eqref{eqn:S2part2} when $q \equiv 1 \imod 4$.\\

Now suppose that $q \equiv 3 \imod 4$ so that $C$ consists of elements in $W$ with partitions $\beta$ such that every $b_i^+$ is odd, and every $b_j^-$ is even. Now each permutation $\tau \in S_{\l}$ with $k$ cycles corresponds to precisely ${2^{\l-k}}$ elements $w$ in $W$ that satisfy the conditions: $\sigma(w)=\tau$,  every odd cycle is positive, and every even cycle is negative. This is exactly the same situation as the case when $q \equiv 1 \imod 4$ and so \eqref{eqn:S2part2} follows in the same way.

The proof  of \eqref{eqn:Sodd} is similar.  If a maximal torus $T$ of \tpe 2 has $k$ direct factors, then by Lemma \ref{lem:QmeetT} we have
\[ \frac{|Q(1) \cap T|}{|T|} = \frac{1}{2^k}.\]
Thus Theorem \ref{thm:maintool} implies that
\begin{displaymath}
\frac{|Q(1)|}{|\GF|}  \ge  \sum_{k=1}^{\l} \frac{c(\l,k)}{4^k\l!}  ,
\end{displaymath}
and \eqref{eqn:Sodd} then follows from \eqref{eqn:clk4^k} and \eqref{eqn:clk4^kasl^3/4}.

This completes the proof for cases where $\GF = \Sp_{2\l}(q)$, or $\SO_{2\l+1}(q)$. If $\GF= \Om_{2n+1}(q)$, then $\GF$ has even index in $\SO_{2\l+1}(q)$, so all odd order elements in $\SO_{2\l+1}(q)$ are contained in $\GF$, and the proportion of odd order elements in $\GF$ is twice the proportion of odd order elements in $\SO_{2\l+1}(q)$. For the proportion of elements of \tpo $2$ in $\GF$, define the quokka set
\[Q':= \{g\in \SO_{2\l+1}(q) \,:\, g \in \Om_{2\l+1}(q) \text{ and } |g|_2=2\}.\]
We take $C$ as before, and we claim that if a maximal torus $T$ of \tpe $2$ has $k$ direct factors, then
\begin{align} \label{eqn:Q'meetT}
\frac{|Q' \cap T|}{|T|}& = \frac{1}{2}- \frac{1}{2^k}.
\end{align}
For $|T \cap \Om_{2\l+1}(q)|= |T|/2$, and the elements in $T$ that are not contained in $\Om_{2\l+1}(q)$ must have even order since $\Om_{2\l+1}(q)$ has index 2 in $\SO_{2\l+1}(q)$. Thus, on the one hand the number of elements in $T$ of \tpo $2$ is $|T|(1- \frac{1}{2^k})$. On the other hand, it is the number of elements in $T$ that are not contained in $\Om_{2\l+1}(q)$ plus the number of elements of \tpo $2$ in $T \cap \Om_{2\l+1}(q)$. So, we have
\begin{displaymath}
|T|\left(1- \frac{1}{2^k}\right) = \frac{|T|}{2} + |Q' \cap T|,
\end{displaymath}
from which \eqref{eqn:Q'meetT} follows easily. Now Theorem \ref{thm:maintool} implies that
\begin{align}
\frac{|Q'|}{|\SO_{2\l+1}(q)|} &\ge \sum_{k=1}^{\l} \frac{c(\l,k)}{2^k\l!}  \left(\frac{2^{k-1}-1}{2^{k}} \right) =   \frac{1}{2}\sum_{k=1}^{\l} \frac{c(\l,k)}{2^k\l!}  - \sum_{k=1}^{\l} \frac{c(\l,k)}{4^k\l!}  \notag\\
& \ge \frac{(2 \l)!}{2^{2\l+1} (\l !)^2} -   \frac{\Gamma(\l+1/4)}{\Gamma(1/4)\l!} \notag
\end{align}
and so the proportion of elements in $\Om_{2\l+1}(q)$ of \tpo $2$ is at least
\begin{align}
\frac{|Q(2)|}{|\Om_{2\l+1}(q)|} &\ge \frac{(2 \l)!}{2^{2\l} (\l !)^2} -   \frac{2\Gamma(\l+1/4)}{\Gamma(1/4)\l!}\ge \frac{25}{29 \sqrt{\pi \l}}-\frac{6}{5 (\l+1)^{3/4}}.
\end{align}
Also, if $\GF = \Sp_{2\l}(q)$ or $\SO_{2\l+1}(q)$ , then the proportion of odd order elements in $\GF/Z(\GF)$ is at least the proportion of odd order elements in $\GF$ by Theorem \ref{thm:projectivequokkasets}, since $Q(1)Z/Z$ is a quokka set and consists of odd order elements in $\GF/Z(\GF)$. Moreover, all of the odd order elements in $\PSO_{2\l+1}(q)$ are contained in $\POm_{2\l+1}(q)$.

Now suppose that $\GF = \PSp_{2\l}(q)$. We define the quokka subset $Q''$ of $\Sp_{2\l}(q)$ to be the set of elements $g$ of \tpo $2$ such that $g^{|g|/2}$ is not central; that is
\begin{displaymath}
Q''= \{ g \in \Sp_{2\l}(q) \,:\, |g|_2=2 \text{ and } g^{|g|/2} \ne -I_{2\l} \}.
\end{displaymath}
Again we use $C$ as before and note that a maximal torus $T$ with \tpe $2$ and $k$ cycles satisfies
\begin{displaymath}
\frac{|Q'' \cap T|}{|T|} \ge 1- \frac{1}{2^k} - \frac{1}{2^k}= 1- \frac{1}{2^{k-1}}
\end{displaymath}
by Lemma \ref{lem:QmeetT}. So Theorem \ref{thm:maintool} implies that
\begin{align} \label{eqn:Spproj}
\frac{|Q''|}{|\Sp_{2\l}(q)|} &\ge \frac{(2 \l)!}{2^{2\l} (\l !)^2} -   \frac{2\Gamma(\l+1/4)}{\Gamma(1/4)\l!}\ge \frac{25}{29 \sqrt{\pi \l}}-\frac{6}{5 (\l+1)^{3/4}}.
\end{align}
Now Theorem \ref{thm:projectivequokkasets} implies that $Q''Z/Z$ is a quokka subset of $\PSp_{2\l}(q)$; it consists of elements of \tpo $2$, and moreover $|Q(2)|/|\GF|$ is at least the  lower bound in \eqref{eqn:Spproj}. Exactly the same argument gives the same lower bound if $\GF=\PSO_{2\l+1}(q)$.

Finally, if $\GF = \POm_{2\l+1}(q)$ then we consider the quokka subset
\begin{displaymath}
Q''':= \{ g \in \SO_{2\l+1}(q) \,:\,g \in \Om_{2\l+1}(q), |g|_2=2 \text{ and } g^{|g|/2} \ne -I_{2\l} \}.
\end{displaymath}
of $\SO_{2\l+1}(q)$. Now for a maximal torus $T \le \SO_{2\l+1}(q)$ of \tpe $2$ with $k$ cycles, we have
\begin{displaymath}
\frac{|Q''' \cap T|}{|T|} \ge \frac{1}{2}- \frac{1}{2^{k-1}}
\end{displaymath}
since the argument when $\GF = \Om_{2\l+1}(q)$ shows that the proportion of elements in $T \cap \Om_{2\l+1}(q)$ that have \tpo $2$ is $1/2-1/2^k$, and also the proportion of elements $g\in T$ with \tpo $2$ such that  $g^{|g|/2} = -I_{2\l}$  is $1/2^k$ by Lemma \ref{lem:QmeetT}(d). It follows from Theorem \ref{thm:maintool} and Theorem \ref{thm:projectivequokkasets} that
\begin{equation} \label{eqn:Oproj}
\frac{|Q(2)|}{|\POm_{2\l+1}(q)|} \ge \frac{(2 \l)!}{2^{2\l} (\l !)^2} -   \frac{4\Gamma(\l+1/4)}{\Gamma(1/4)\l!}\ge \frac{25}{29 \sqrt{\pi \l}}-\frac{12}{5 (\l+1)^{3/4}}. \qedhere
\end{equation}
\end{proof}

\begin{pro} \label{pro:caseS2parts2j}
Let $\GF = \mathrm{Sp}_{2\l}(q)$ or $\SO_{2\l+1}(q)$, and define $t$ as in \eqref{eqn:introdeft}. If $1 \le 2^{\a}/t \le n$, then the proportion of elements that have \tpo $2^{\a}$ is
\begin{align} \label{S1}
\frac{|Q(2^{\a})|}{|\GF|} \ge \frac{p_{\l}(2^{\a}/t)}{4}
\end{align}
where $p_{\l}(2^{\j})$ is the proportion of elements in $S_{\l}$ of \tpo $2^{\j}$. If $2^b/t$ does not divide $n$, then \eqref{S1} holds for $\GF/Z(\GF)$.
\end{pro}
\begin{proof}
We use Theorem \ref{thm:maintool} with $C \subset W$ consisting of all $F$-classes that correspond to maximal tori with \tpe $2^{\a}$. By Lemma \ref{caseS2parts},  if $2 \le 2^{\a}/t \le n$, then at least half of the elements $w$ in $W$, such that $\sigma(w)$ in $S_{\l}$ has \tpo $2^{\a}/t$, have a positive cycle of \tpo $2^{\a}/t$, and therefore correspond to maximal tori $T$ with \tpe $2^{\a}$. Furthermore, by Lemma \ref{lem:QmeetT}, at least half of the elements in a maximal torus $T$ of \tpe $2^{\a}$ have \tpo $2^{\a}$. Thus, Theorem \ref{thm:maintool} implies that
\begin{displaymath}
\frac{|Q(2^{\a})|}{|\GF|} \ge \frac{p_{\l}(2^{\a}/t)}{4}.
\end{displaymath}

If $q \equiv 1 \imod 4$ and $2^{\a}=t$, then maximal tori of $\GF$ of \tpe $t=(q-1)_2$ correspond to partitions $\beta$, where all even cycles are negative, and there is at least one positive (odd) cycle. If $q \equiv 3 \imod 4$ then $C$ contains elements in $W$, where all of the even cycles are negative, and there is at least one negative odd cycle. In both cases, at least half of the elements in $W$ that correspond to a permutation in $S_{\l}$ of odd order satisfy this condition. The proportion of such elements in $W$ is therefore at least $\frac{p_{\l}(1)}{2}$.
Now using Lemma \ref{lem:QmeetT} and Theorem \ref{thm:maintool} we have
\begin{displaymath}
\frac{|Q(t)|}{|\GF|}\ge \frac{p_{\l}(1)}{4}.
\end{displaymath}
The same argument as for Proposition \ref{pro:gln} shows that \eqref{S1} holds for $G/Z(G)$ when $2^{\j}/t$ does not divide $\l$.
\end{proof}
\begin{remark} \label{rem:6.4}
\emph{
A similar argument to the one in Remark \ref{lin:semiregular} yields that if $2^\a/t \ge 2$ and $2^\a/t$ divides $\l$, then
\begin{displaymath}
\frac{|Q(2^\a)|}{|\GF|}\ge \frac{p_{\l}(1)}{4} - \frac{3}{4\l^{1-\frac{t}{2^{b+1}}}}
\end{displaymath}
for $G=\PSp_{2\l}(q)$ or  $\PSO_{2\l+1}(q)$. This completes the proof of Theorem \ref{thm:constantclassical} in the symplectic and odd-dimensional orthogonal cases.}
\end{remark}
\section{Orthogonal case}
The structure of the maximal tori in $\SO^{\epsilon}_{2\l}(q)$ is essentially the same as in the symplectic case, but the Weyl group has index 2 in $W_{B}:=C_2 \wr S_{\l}$. Moreover, in the untwisted case, an $F$-conjugacy class in $W$ is a conjugacy class of permutations in $W_B$ with an even number of negative cycles. In the twisted case, an $F$-conjugacy class $C$ in $W$ corresponds to a conjugacy class $C^{\prime}$ of permutations in $W_B$ with an odd number of negative cycles; if we define the $2$-cycle $w=(\l \, \l^{\prime}) \in W_B$, then $C=wC^{\prime}$.

 In Proposition \ref{pro:caseO2}, we consider proportions of elements in even dimensional orthogonal groups with \tpo $1$ or $2$, which completes the proof of Theorem \ref{thm:oddandtwiceodd}.

\begin{pro} \label{pro:caseO2}
Let $(\GF,\delta_1,\delta_2) = (\SO_{2\l}^{\epsilon}(q),1,1)$, $(\Om_{2\l}^{\epsilon}(q),2,2)$, $(\PSO_{2\l}(q),1,2)$, or $(\POm^{\epsilon}_{2\l}(q),2,4)$. If $Q(1):= \{ g \in \GF \,:\, |g|_2=1\}$, then the proportion of odd order elements in $\GF$ is
\begin{align} \label{OQ1}
\frac{|Q(1)|}{|\GF|} \ge \left\{
  \begin{array}{ll}
   \frac{\Gamma(\l-1/4)\delta_1}{2\Gamma(3/4)}+\frac{\Gamma(\l+1/4)\delta_1}{2\Gamma(1/4)}, & \hbox{if $q \equiv 1 \imod 4$, and $\epsilon=-$; or} \\
& \hbox{$q \equiv 3 \imod 4$, $\epsilon=+$, and $\l$ odd; or} \\
 & \hbox{$q \equiv 3 \imod 4$, $\epsilon=-$, and $\l$ even} \\
\frac{\Gamma(\l+1/4)\delta_1}{2\Gamma(1/4)}-\frac{\Gamma(\l-1/4)\delta_1}{2\Gamma(3/4)}, & \hbox{if $q \equiv 1 \imod 4$ and $\epsilon=+$; or}\\
 & \hbox{$q \equiv 3 \imod 4$, $\epsilon=+$, and $\l$ even; or} \\
 & \hbox{$q \equiv 3 \imod 4$, $\epsilon=-$, and $\l$ odd.}
  \end{array}
\right.
\end{align}
Simplified lower bounds for these estimates are $\frac{\delta_1}{4(\l+1)^{3/4}}+\frac{\delta_1}{5(\l+1)^{5/4}}$ and $\frac{\delta_1}{4(\l+1)^{3/4}}-\frac{18\delta_1}{25(\l+1)^{5/4}}$ respectively.
Also, if $Q(2):= \{ g \in \GF \,:\, |g|_2=2\}$, then the proportion of elements in $\GF$ of \tpo $2$ is
\begin{align} \label{OQ2}
\!\!\frac{|Q(2)|}{|\GF|} \! \ge \! \left\{
  \begin{array}{ll}
    \! \frac{(2\l)!}{2^{2\l}(\l !)^2} \frac{2\l}{2\l-1} \! - \! \frac{\Gamma(\l \! -\! \frac{1}{4})\delta_2}{2\Gamma(3/4)} \!-\! \frac{\Gamma(\l\!+\! \frac{1}{4})\delta_2}{2\Gamma(1/4)}, & \!\!\!\hbox{if $q \equiv 1 \!\! \imod 4$, and $\epsilon\!=\!-$; or} \\
& \!\!\! \hbox{$q \equiv 3 \!\!\imod 4$, $\epsilon\!=\!+$, and $\l$ odd; or} \\
 & \!\!\! \hbox{$q \equiv 3 \!\! \imod 4$, $\epsilon\!=\!-$, and $\l$ even} \\
\! \frac{(2\l)!}{2^{2\l}(\l !)^2} \frac{2\l \! -\! 2}{2\l \!-\!1} \!+ \! \frac{\Gamma(\l \!- \! \frac{1}{4})\delta_2}{2\Gamma(3/4)}\! -\! \frac{\Gamma(\l \! + \!\frac{1}{4})\delta_2}{2\Gamma(1/4)}, & \!\! \hbox{if $q \equiv 1 \!\! \imod 4$, and $\epsilon \!= \!+$; or }\\
 & \!\! \hbox{$q \equiv 3 \!\!\imod 4$, $\epsilon \!=\!+$, and $\l$ even; or} \\
 &\!\! \hbox{$q \equiv 3 \!\! \imod 4$, $\epsilon \!=\!-$, and $\l$ odd.}
  \end{array}
\right.
\end{align}
Simplified lower bounds for these estimates are $\frac{50\l}{29(2\l-1) \sqrt{\pi\l}}-\frac{3\delta_2}{10(\l+1)^{3/4}} - \frac{9\delta_2}{25(\l+1)^{5/4}}$ and $\frac{50(\l-1)}{29(2\l-1) \sqrt{\pi\l}}-\frac{3\delta_2}{10(\l+1)^{3/4}} + \frac{\delta_2}{10(\l+1)^{5/4}}$ respectively.
\end{pro}
\begin{proof}
First of all, we suppose that $G= \SO^\epsilon_{2\l}(q)$. We prove both inequalities using Theorem \ref{thm:maintool}, where $C \subset W$ is the union of $F$-classes that correspond to maximal tori of \tpe $2$.
  Suppose that $q \equiv 1 \imod 4$ and $\epsilon=+$. Then $C$ consists of $w \in W$ such that all cycles of  $w \in W$ are negative cycles. Since the number of negative cycles must be even, $\sigma(w) \in S_{\l}$ must have an even number of cycles. The proportion of elements in $w\in W$ with $k$ negative cycles and no positive cycles, where $k$ is even, is
\begin{displaymath}
\frac{c(\l,k)2^{\l-k}}{\l! 2^{\l-1}} = \frac{2c(\l,k)}{\l ! 2^k }.
\end{displaymath}
Therefore, by Theorem \ref{thm:maintool}, we have
\begin{displaymath}
\frac{|Q(2)|}{|\GF|} \ge \sum_{\substack{{k=1}\\{k \text{ even}}}}^{\l} \frac{2c(\l,k)}{\l ! 2^k} \left(1-\frac{1}{2^{k}} \right) = 2\sum_{\substack{{k=1}\\{k \text{ even}}}}^{\l} \frac{c(\l,k)}{\l ! 2^k}- 2\sum_{\substack{{k=1}\\{k \text{ even}}}}^{\l} \frac{c(\l,k)}{\l ! 4^k}
\end{displaymath}
and
\begin{displaymath}
\frac{|Q(1)|}{|\GF|} \ge \sum_{\substack{{k=1}\\{k \text{ even}}}}^{\l} \frac{2c(\l,k)}{\l ! 2^k} \left(\frac{1}{2^{k}} \right) = 2\sum_{\substack{{k=1}\\{k \text{ even}}}}^{\l} \frac{c(\l,k)}{\l ! 4^k}
\end{displaymath}
for $q \equiv 1 \imod 4$ and $\epsilon=+$. We calculated these sums in \eqref{eqn:stirling2^keven} and \eqref{eqn:stirling4^keven}, and we can find the simplified lower bounds using \eqref{eqn:2ellfactorial} and \eqref{eqn:boundsclk4^keven}.\\
The other cases are similar. Suppose that $q \equiv 1 \imod 4$ and $\epsilon=-$. Now $C$ consists of $w \in W$ such that all of the cycles of $w$ are negative. Since the number of negative cycles must be odd when $\epsilon=-$, $\sigma(w) \in S_{\l}$ must have an odd number of cycles. Moreover, by Theorem \ref{thm:maintool} we have
\begin{align*}
\frac{|Q(2)|}{|\GF|} &\ge  \sum_{\substack{{k=1}\\{k \text{ odd}}}}^{\l} \frac{2c(\l,k)}{\l ! 2^k}\left(1- \frac{1}{2^k}\right)\\
&\ge 2\sum_{\substack{{k=1}\\{k \text{ odd}}}}^{\l} \frac{c(\l,k)}{\l ! 2^k}   - 2\sum_{\substack{{k=1}\\{k \text{ odd}}}}^{\l} \frac{c(\l,k)}{\l ! 4^k}.
\end{align*}
Now we can use \eqref{eqn:stirling2^kodd} and \eqref{eqn:stirling4^kodd} to obtain the required lower bounds. We obtain the simplified lower bounds using \eqref{eqn:2ellfactorial} and \eqref{eqn:boundsclk4^kodd}. which completes the proof for the case where $q \equiv 1 \imod 4$ and $\epsilon=-$.

Next, suppose that $q \equiv 3 \imod 4$ and $\epsilon=+$. Then $C$ consists of elements $w \in W$ such that all odd length cycles of $w$ are positive, and all even length cycles of $w$ are negative. Since $\epsilon=+$, the number of negative (even length) cycles must be even, and thus $\sigma(w) \in S_{\l}$ must have an even number of cycles of even length. Recall that we defined $c_1(\l,k)$ to be the number of permutations in $S_{\l}$ with $k$ cycles, and an even number of even length cycles; $c_2(\l,k)$ was defined to be the number of permutations in $S_{\l}$ with $k$ cycles, and an odd number of even length cycles.  By Theorem \ref{thm:maintool}, we have
\begin{align*}
\frac{|Q(2)|}{|\GF|} &\ge 2\sum_{k=1}^{\l} \frac{c_1(\l,k)}{\l ! 2^k}- 2 \sum_{k=1}^{\l} \frac{c_1(\l,k)}{\l ! 4^k}
\end{align*}
and
\begin{align*}
\frac{|Q(1)|}{|\GF|} &\ge 2\sum_{k=1}^{\l} \frac{c_1(\l,k)}{\l ! 4^k}.
\end{align*}
Thus if $\l$ is even, $q \equiv 3 \imod 4$ and $\epsilon=+$, then Lemma \ref{lem:numberevencycles} implies that
\begin{align*}
\frac{|Q(2)|}{|\GF|}&\ge 2 \sum_{\substack{{k=1}\\{k \text{ even}}}}^{\l} \frac{c(\l,k)}{\l ! 2^k}- 2\sum_{\substack{{k=1}\\{k \text{ even}}}}^{\l} \frac{c(\l,k)}{\l ! 4^k},
\end{align*}
and
\begin{align*}
\frac{|Q(1)|}{|\GF|}&\ge 2\sum_{\substack{{k=1}\\{k \text{ even}}}}^{\l} \frac{c(\l,k)}{\l ! 4^k}.
\end{align*}
We can now obtain the required lower bounds using  \eqref{eqn:stirling2^keven}, \eqref{eqn:stirling4^keven}, \eqref{eqn:2ellfactorial}, and \eqref{eqn:boundsclk4^keven}.
If $\l$ is odd, $q \equiv 3 \imod 4$ and $\epsilon=+$, then
\begin{align*}
\frac{|Q(2)|}{|\GF|}&\ge 2\sum_{\substack{{k=1}\\{k \text{ odd}}}}^{\l} \frac{c(\l,k)}{\l ! 2^k}- 2\sum_{\substack{{k=1}\\{k \text{ odd}}}}^{\l} \frac{c(\l,k)}{\l ! 4^k},
\end{align*}
and
\begin{align*}
\frac{|Q(1)|}{|\GF|}&\ge 2\sum_{\substack{{k=1}\\{k \text{ odd}}}}^{\l} \frac{c(\l,k)}{\l ! 4^k}.
\end{align*}
Now we obtain the lower bounds as before using \eqref{eqn:stirling2^kodd}, \eqref{eqn:stirling4^kodd}, \eqref{eqn:2ellfactorial}, and \eqref{eqn:boundsclk4^kodd}.
Finally, suppose that $q \equiv 3 \imod 4$ and $\epsilon=-$. Then  $C$ consists of $w \in W$ such that all of the odd length cycles of $w$ are positive and all of the even length cycles of $w$ are negative. Since $\epsilon=-$, the number of negative (even length) cycles must be odd. Using Lemma \ref{lem:numberevencycles} and Theorem \ref{thm:maintool}, we have
\begin{align*}
\frac{|Q(2)|}{|\GF|}&\ge \sum_{k=1}^{\l} \frac{2c_2(\l,k)}{\l ! 2^k}\left(1- \frac{1}{2^k}\right);
\end{align*}
if $\l$ is even, then
\begin{align*}
\frac{|Q(2)|}{|\GF|}&\ge 2\sum_{\substack{{k=1}\\{k \text{ odd}}}}^{\l} \frac{c(\l,k)}{\l ! 2^k}- 2\sum_{\substack{{k=1}\\{k \text{ odd}}}}^{\l} \frac{c(\l,k)}{\l ! 4^k},
\end{align*}
and
\begin{align*}
\frac{|Q(1)|}{|\GF|}&\ge 2\sum_{\substack{{k=1}\\{k \text{ odd}}}}^{\l} \frac{c(\l,k)}{\l ! 4^k};
\end{align*}
if $\l$ is odd then
\begin{align*}
\frac{|Q(2)|}{|\GF|}&\ge 2\sum_{\substack{{k=1}\\{k \text{ even}}}}^{\l} \frac{c(\l,k)}{\l ! 2^k}- 2\sum_{\substack{{k=1}\\{k \text{ even}}}}^{\l} \frac{c(\l,k)}{\l ! 4^k},
\end{align*}
and
\begin{align*}
\frac{|Q(1)|}{|\GF|}&\ge 2\sum_{\substack{{k=1}\\{k \text{ even}}}}^{\l} \frac{c(\l,k)}{\l ! 4^k}.
\end{align*}
We proceed as in the other cases. This completes the proof of Proposition \ref{pro:caseO2} for the case $q \equiv 3 \imod 4$ and $\epsilon=-$, and this completes the proof for $\GF=\SO^\epsilon_{2\l}(q)$.

If $\GF= \Om_{2n}^\epsilon(q)$, then $\GF$ has even index in $\SO_{2n}(q)$, so all odd order elements in $\SO_{2n}^\epsilon(q)$ are contained in $\GF$, and the proportion of odd order elements in $\GF$ is twice the proportion of odd order elements in $\SO_{2n}^\epsilon(q)$. Also, the proportion of odd order elements in $\GF/Z(\GF)$ is at least the proportion of odd order elements in $\GF$ by Theorem \ref{thm:projectivequokkasets}, since $Q(1)Z/Z$ is a quokka set and consists of odd order elements in $\GF/Z(\GF)$. Thus \eqref{OQ1} is proved in all cases.

For\eqref{OQ2}, we adapt our proof from the case $\GF=\SO^\epsilon_{2\l}(q)$ as we did in Proposition \ref{caseS2part2}. If $\GF= \Om_{2n}^\epsilon(q)$ then we let
\[Q':= \{g\in \SO^{\epsilon}_{2\l}(q) \,:\, g \in \Om^{\epsilon}_{2\l}(q) \text{ and } |g|_2=2\};\]
if $\GF = \PSO^\epsilon_{2\l}(q)$, then we let
\begin{displaymath}
Q''= \{ g \in \SO^{\epsilon}_{2\l}(q) \,:\, |g|_2=2 \text{ and } g^{|g|/2} \ne -I_{2\l} \};
\end{displaymath}
if $\GF = \POm_{2\l}^{\epsilon}(q)$ then we let
\begin{displaymath}
Q''':= \{ g \in \SO_{2\l}(q) \,:\,g \in \Om^{\epsilon}_{2\l}(q), |g|_2=2 \text{ and } g^{|g|/2} \ne -I_{2\l} \}.
\end{displaymath}
We note that if $T$ is a maximal torus in $\SO_{2\l}^{\epsilon}$ of \tpe $2$ corresponding to a partition with $k$ parts, then
\begin{align*}
\frac{|Q' \cap T|}{|T|}& = \frac{1}{2}- \frac{1}{2^k}, \\
\end{align*}
and by Lemma \ref{lem:QmeetT}(d),
\begin{align*}
\frac{|Q'' \cap T|}{|T|} &\ge 1- \frac{1}{2^{k-1}},
\end{align*}
and
\begin{align*}
\frac{|Q''' \cap T|}{|T|} &\ge \frac{1}{2} - \frac{1}{2^{k-1}}.
\end{align*}
We use Theorem \ref{thm:maintool} with $C$ as before and we apply Theorem \ref{thm:projectivequokkasets} for $Q''$ and $Q'''$, noting that if $Z=Z(\GF)$ in these cases, then  $Q''Z/Z$, and  $Q'''Z/Z$ consist of elements of \tpo $2$ in $G/Z$. Now \eqref{OQ2} follows for all cases.
\end{proof}

\begin{pro} \label{pro:caseO2jp}
Let $\GF = \SO^{\epsilon}_{2\l}(q)$, where $\l \ge 2$, and define $t$ as in \eqref{eqn:introdeft}. If $1 \le 2^{\a}/t < n$ and $\l_2 \ne 2^{\a}/t$, then the proportion of elements in $\GF$ with \tpo equal to $2^{\a}$ is
\begin{align}\label{O1}
\frac{|Q(2^{\a})|}{|\GF|} \ge \frac{p_{\l}(2^{\a}/t)}{4},
\end{align}
where $p_{\l}(2^{\j}/t)$ is the proportion of elements in $S_{\l}$ of \tpo $2^{\j}/t$. If $\l_2=2^{\a}/t$, then
the proportion of elements in $\GF$ with \tpo equal to $2^{\a}$ is
\begin{align}\label{O1:2tobdividesn}
\frac{|Q(2^{\a})|}{|\GF|} \ge \frac{p_{\l}(2^{\a}/t)}{4} - \frac{1}{4\l}.
\end{align}
If $2^{\a}/t$ does not divide $\l$, then \eqref{O1} holds for $G=\PSO^{\epsilon}_{2\l}(q)$ as well.
\end{pro}
\begin{proof}
As in the symplectic case, we use Theorem \ref{thm:maintool}, where $C$ is the union of $F$-classes that correspond to maximal tori  with \tpe  $2^{\a}$. First suppose that $2 \le 2^{\a}/t < n$ and $n_2 \ne 2^{\a}/t$. Suppose that $\sigma(w) \in S_{\l}$ has \tpo  $2^{\a}/t$. Then at least half of the elements in $W_B$ with this cycle structure in $S_{\l}$ have a positive cycle of length $2^{\a}d/t$, where $d$ is an odd integer. Since $\l_2 \ne 2^\a/t$, $\sigma(w)$ must have at least one other cycle, and we claim that precisely half of these elements in $W_B$ have an even number of negative cycles and thus correspond to $F$-classes of $W$ in the $\epsilon=+$ case; the other half have an odd number of negative cycles and correspond to $F$-classes of $W$ in the $\epsilon=-$ case. For suppose that $\sigma(w)$ has $m+1$ cycles; without loss of generality, the first one has length $2^{\j}d$, and this has been assigned a positive sign in $W$. For each of the remaining $m$ cycles, precisely half of the $2^{\l}$ corresponding elements in $W_B$ have that cycle assigned as positive. Furthermore, there are $2^m$ possible ways to distribute positive and negative signs to the $m$ cycles. Of these $2^m$ ways, there are ${m \choose i}$ ways that result in $i$ negative signs. And so the number of ways resulting in an even number of negative cycles is
\begin{displaymath}
\sum_{\substack{{i=0}\\{i \text{ even}}}}^m {m \choose i} =2^{m-1},
\end{displaymath}
which proves the claim. Therefore in both the $\epsilon=+$ and $\epsilon=-$ cases, at least a quarter of the elements in $w \in W_B$ with $\sigma(w)$ of \tpo $2^{\a}/t$ (where $n_2 \ne 2^{\a}/t$), correspond to maximal tori $T$ with \tpe $2^{\a}$. By Lemma \ref{lem:QmeetT}, at least half of the elements in $T$ have \tpo $2^{\a}$. Thus by Theorem \ref{thm:maintool}, we have
\begin{displaymath}
\frac{|Q(2^{\a})|}{|\GF|} \ge \frac{p_{\l}(2^{\a}/t)}{4}.
\end{displaymath}
If $2^{\a}/t$ does not divide $\l$ (in particular $\l_2 \ne 2^{\a}/t$), then the same argument as for Proposition \ref{pro:gln} shows that \eqref{O1} holds for $\PSO_{2\l}^{\epsilon}(q)$ as well.
\par
If $\l_2=2^\a/t$, then some of the permutations $\sigma(w)$ in $S_\l$ of \tpo $2^\a/t$ will be $\l$-cycles, for which the argument above fails. So we repeat the argument above but exclude those $w \in W$ for which $\sigma(w) \in S_\l$ is an $\l$-cycle. Since the proportion of $\l$-cycles in $S_\l$ is $\frac{1}{\l}$,  Theorem \ref{thm:maintool} yields
\begin{displaymath}
\frac{|Q(2^{\a})|}{|\GF|} \ge \frac{p_{\l}(2^{\a}/t)}{4}- \frac{1}{4\l}.
\end{displaymath}
If $2^{\a}=t$ and $q \equiv 1 \imod 4$, then $C$ consists of elements in $W$ with at least one positive odd cycle, and with all even cycles being negative. If $q \equiv 3 \imod 4$, then the condition is that there is at least one negative cycle of odd length and all of the even length cycles are negative. As in the previous case, given a permutation $\tau$ in $S_{\l}$ that is not an $\l$-cycle, precisely half of the elements $w \in W_B$ such that $\sigma(w)=\tau$ have an even number of negative cycles. Thus for both $\epsilon=+$ and $\epsilon=-$, if $\l$ is even, then the proportion of elements in $W$ that are contained in $C$  is at least $\frac{p_{\l}(1)}{2}$ 
 and Lemma \ref{lem:QmeetT} and Theorem \ref{thm:maintool} yield
\begin{displaymath}
\frac{|Q(t)|}{|\GF|} \ge  \frac{p_{\l}(1)}{4}.
\end{displaymath}
If $n$ is odd (so that $n_2=2^\a/t=1$), then we argue as above, but exclude those $w \in W$ for which $\sigma(w)$ is an $\l$-cycle, and we have
\begin{equation}
\displaystyle \frac{|Q(t)|}{|\GF|} \ge \displaystyle  \frac{p_{\l}(1)}{4} - \frac{1}{4\l}. \qedhere
\end{equation}
\end{proof}
\begin{remark} \label{rem:7.3}
\emph{
A similar argument to the one in Remark \ref{lin:semiregular} yields that if $2^\a/t \ge 2$,  $2^\a/t$ properly divides $\l$, and $Q(2^b)$ denotes the proportion of elements of \tpo $2^b$ in $\PSO_{2\l}^{\epsilon}(q)$, then
\begin{displaymath}
\frac{|Q(2^\a)|}{|\PSO_{2\l}^{\epsilon}(q)|}\ge \frac{p_{\l}(1)}{4} - \frac{3}{4\l^{1-\frac{t}{2^{b+1}}}}.
\end{displaymath}
This completes the proof of Theorem \ref{thm:constantclassical}.}
\end{remark}

\section{Strong involutions}
First, we recall the definition of a strong involution.
\begin{dfn}
\emph{
Let $z$ be an involution in a classical group $G$, with natural module of dimension $d$.  We say that $z$ is a \textit{strong involution} if the dimension of the fixed point space of $z$ is contained in the interval $[d/3,2d/3)$.}
\end{dfn}
Consider an involution $z$ in $G$ with centralizer $C = Cl_m(q)\! \times \! Cl_{d-m}(q)$, where $Cl_m(q)$ and $Cl_{d-m}(q)$ are classical groups of dimension $m$ and $d-m$. If $(x,y)$ is chosen at random in $C$ then let
\begin{equation} \label{e:PMG}
  P_m^G := \mathbb{P}(r:=|(x,y)| \text{ is even and } (x,y)^{r/2} \in Cl_m(q) \times \{ 1 \}).
\end{equation}
We will often write $P_m$ for $P_m^G$ when there is little possibility of confusion.

The following result is of interest for the problem of constructing the direct factors in $C$.
\begin{pro} \label{pro:balanced}
Let $\GF$ be one of the groups in column 1 of Table \ref{tab:forbalanced} and suppose $m$ and $d$ are positive integers satisfying the conditions in column $2$.   Then a lower bound for $P_m^G$ is given in column $3$ of Table \ref{tab:forbalanced}.
\end{pro}
\begin{table}[htdp]
\begin{center}
 \caption{Lower bounds for $P_m^G$ in Proposition \ref{pro:balanced}}
\begin{tabular}{ccc}
\hline
 \vspace*{-3.5mm}\\
 $G$ & Conditions & Lower bound on $P_m^G$ \\
 \vspace*{-4mm}\\
 \hline
 \vspace*{-3.5mm}\\
 $\GL^{\epsilon}_d(q)$& $\frac{d}{3} \le m \le \frac{d}{2}$, $m \ge 2$ & $0.078125$ \\
 $\Sp_d(q)$ & $\frac{d}{3} \le m \le \frac{d}{2}$, $m \ge 2$ & $0.039062$  \\
  $\SO_d(q)$ & $\frac{d}{3} \le m \le \frac{d}{2}$, $m \ge 3$ & $0.039062$  \\
  $\GL^{\epsilon}_d(q)$& $\frac{d}{2} \le m \le \frac{2d}{3}$ & $0.087683$ \\ 
 $\Sp_d(q)$, $\SO^{\epsilon}_d(q)$  & $\frac{d}{2} \le m \le \frac{2d}{3}$ & $0.043216$  \\
\vspace*{-3.5mm}\\
\hline
\vspace*{-3.5mm}\\
 \end{tabular} 
\end{center}
\label{tab:forbalanced}
\end{table}

\begin{remark}
\emph{(i) 
A MAGMA calculation shows that the proportion above in $\GL_3(3) \times \GL_6(3)$ is 0.1857852; so our lower bound fails to be sharp by a factor of at most 2.38 in the linear case. For the case where $d/2 \le m \le 2d/3$ we find that the proportion in $\GL_2(3) \times \GL_2(3)$ is 0.25781250.}\\
\emph{(ii)
The choice of $d/3$ and $d/2$ as lower and upper bounds for $m$ can be changed; essentially to any $C d$ and $D d$ where $0< C < D <1$, and we would obtain alternative lower bounds (here we say that $z$ is a $(C,D)$-balanced involution). For example, if we require that  $d/6 \le m \le d/2$  (taking $z$ to be a so-called semistrong involution), then our lower bounds\footnote{corresponding to $d=18$, and $m=3$} are 0.052368 when $G=\GL_d^{\epsilon}(q)$
 and 0.0276370 when $G\ne \GL_d^{\epsilon}(q)$; if we require that  $d/2 \le m \le 5d/6$, then we obtain lower bounds of $0.0870225050$
and $0.0430945858$ respectively.
}
\end{remark}
To help to prove Proposition \ref{pro:balanced}, we find lower bounds on the proportion of elements in $Cl_{d-m}(q)$ with \tpo strictly less than $2^b$.
\begin{lemma} \label{lem:notdivisibleby}
Let $\GF = \GL^{\epsilon}_{n}(q)$, $\Sp_{2n}(q)$, $\SO_{2n+1}(q)$, or $\SO^{\epsilon}_{2n}(q)$, where $q$ is odd, and define $t$ as in \eqref{eqn:introdeft}. If $2^\a \ge 2t$, then the proportion of elements in $\GF$ of order not divisible by $2^\a$ is at least the proportion $s_{\neg 2^\a/t}(n)$ of elements in $S_{n}$ of order not divisible by $2^{\a}/t$.
\end{lemma}
\begin{proof}
Let $Q$ be the subset of $\GF$ consisting of all elements in $\GF$ with order not divisible by $2^b$. Then clearly $Q$ is a quokka set. We use Theorem \ref{thm:maintool} where $C$ is the union of $F$-classes in $W$ that correspond to maximal tori with \tpes strictly less than $2^\a$. So in the linear and unitary cases, $C$ consists of $w \in W=S_n$ of order not divisible by $2^\a/t$. In the other cases, $C$ consists of $w \in W$ such that $\sigma(w) \in S_n$ has order  not divisible by $2^\a/t$.
In all cases, we can apply Theorem \ref{thm:maintool} with $A=1$, and the result follows easily.
\end{proof}

\begin{proof}[Proof of Proposition \ref{pro:balanced}]
First suppose that $\GF = \GL^{\epsilon}_n(q)$ so that $d=n$ and  $C = \GL_m(q) \times \GL_{n-m}(q)$. Now an element in $C$ is of the form $(x,y)$ with $x \in \GL_m(q)$ and $y \in \GL_{n-m}(q)$. We seek a lower bound on the proportion $P_m$ of elements $(x,y) \in C$ such that $|x|_2> |y|_2$.  For $ m\le 399$, ($4\le n \le  1200$) a MAGMA calculation shows that the proportion of elements in
 $S_m \times S_{n-m}$, where $n/3 \le m \le n/2$,  with \tpo in $S_m$ greater than in $S_{n-m}$, is at
 least
\begin{displaymath}
 \sum_{2 \le 2^{\j} \le m} p_m(2^{\j}) s_{\neg 2^{\a}}(n-m) \ge 0.156250.
\end{displaymath}
Thus, by Proposition \ref{pro:gln} and Lemma \ref{lem:notdivisibleby}, if $m \le 399$, then
\begin{displaymath}
P_m \ge \sum_{2 \le 2^{\j} \le m} \frac{p_m(2^{\j})}{2}  s_{\neg 2^{\a}}(n-m) \ge  \frac{0.156250}{2} = 0.078125.
\end{displaymath}
  So let us assume that $m \ge 400$ (and $n \ge 800$) and we can use the lower bounds in Table \ref{tab:forpropenot2e}. Let $\I$ be the collection of intervals \[\I:=\left\{[c,2c) \, : \, c \in \{ a/4,a/2,a,2a,4a,8a,16a\}\right\}.\]
Now the proportion $P_m$ of elements $(x,y) \in C$ such that $|x|_2> |y|_2$ is at least the sum over $\I$ of the proportion of  elements $(x,y) \in C$ such that
$|x|_2$ is the unique $2$-power $2^\a=\alpha_c \log(m)t$ with $\alpha_c \in [c,2c)$, and   $\alpha_c \log(m)t$ does not divide $|y|$.

By Lemma \ref{lem:notdivisibleby} and \cite{BLNPS}, the proportion of elements $y$ in $Cl_{n-m}(q)$ not divisible by $\alpha_c \log(m)t$, where $\alpha_c \log(m)\ge 2$, is at least
\begin{align*}
s_{\neg \alpha_c \log(m)}(n-m) &\ge c(\alpha_c \log(m))(1-1/(n-m))(n-m)^{- \frac{1}{\alpha_c \log(m)}}\\
&\ge c(\alpha_c \log(m))(1-1/m)(2m)^{- \frac{1}{\alpha_c \log(m)}} \\
&\ge c(\alpha_c \log(m))(399/400)e^{-\frac{1}{\alpha_c}- \frac{\log(2)}{\alpha_c \log(m)}}.
\end{align*}
In other words, in the linear and unitary cases,  if $m \ge 400$, then
\begin{align*}
P_m &\ge \sum_{[c,2c) \in \I} s_{\neg \alpha_c \log(m)}(n-m) p_n(\alpha_c \log(m))/2 \\
&\ge \frac{1}{2}  \sum_{[c,2c) \in \I} \frac{399c(\alpha_c \log(m))}{400}e^{-\frac{1}{\alpha_c}\left(1+ \frac{\log(2)}{\log(400)}\right)} p_n(\alpha_c \log(m)).
\end{align*}
We record in Table \ref{tab:balanced} the calculation  of a lower bound for the sum above, where we use the lower bounds on $p_n(\alpha_c\log(m))$ from Table \ref{tab:forpropenot2e}.
So when we are working with $\GL_m(q)\! \times \!\GL_{n-m}(q) \le \GL_n(q)$ or $\GU_m(q) \!\times \!\GU_{n-m}(q) \le \GU_n(q)$ the proportion $P_m$, for $m \ge 400$,  is bounded below by
\begin{align*}
P_m& \ge& \sum_{[c,2c) \in \I} s_{\neg \alpha_c \log(m)}(n-m) p_n(\alpha_c \log(m))/2   \ge  \frac{0.1592612385}{4}
 = 0.0796306192.
\end{align*}
Therefore, for all $m$ in the linear and unitary cases
$P_m \ge 0.078125$,
thus proving Proposition \ref{pro:balanced} in these cases.

\begin{table}[htbp]
\caption{Strong involution calculation, $m \ge 400$.} \label{tab:balanced}
\vspace*{+1mm}
\begin{tabular}{cccc}
\hline
\vspace*{-4mm}\\
\multicolumn{1}{c}{Interval $I$}  & \multicolumn{1}{c}{Lower bound on} & \multicolumn{1}{c}{Lower bound } & \multicolumn{1}{c}{Contribution} \\
\multicolumn{1}{c}{containing $\alpha$}  & \multicolumn{1}{c}{$e^{-\frac{1}{\alpha(1+\log(2)/\log(400))}}$} & \multicolumn{1}{c}{from Table \ref{tab:forpropenot2e}} & \multicolumn{1}{c}{} \\
\vspace*{-4mm}\\
\hline
\vspace*{-4mm}\\
$[a/4,a/2)$ & $0.0136357218$ & 0.1170040878 & 0.0012697907
 \\
$[a/2,a)$ & 0.1167720933 & 0.2351203791 & 0.0273868601 \\
$[a,2a)$ & 0.3417193194 & 0.1531015975 & 0.0521869793 \\
$[2a,4a)$ & 0.5845676346 & 0.0605468750 & 0.0353052591 \\
$[4a,8a)$ & 0.7645702287 & 0.0307617188 & 0.0234606956 \\
$[8a,16a)$ & 0.8743970658 & 0.0155029297 & 0.0135218269 \\
$[16a,32a)$ & 0.9350920093 & 0.0065085753 & 0.0060709014 \\ \hline
\vspace*{-4mm}\\
 &    &  & 0.1592023132  \\
 \vspace*{-4mm}\\
 \cline{4-4}
\end{tabular}
\end{table}
When we are working with $\mathrm{Sp}_m(q) \!\times \!\mathrm{Sp}_{d-m}(q) \!\le \mathrm{Sp}_d(q)$, we find a lower bound for the proportion $P_m$ in the same way. Here  $d = 2\l$ and we write $m = 2k$. Using Proposition \ref{pro:caseS2parts2j}, if $k \ge 400$, then
\begin{align*}
P_m &\ge \sum_{[c,2c) \in \I} s_{\neg \alpha_c \log(k)}(\l-k) p_{k}(\alpha_c \log(k))/4   \\
&\ge 1/4  \sum_{[c,2c) \in \I} \frac{399c(\alpha_c \log(k))}{400}e^{-\frac{1}{\alpha_c}\left(1+ \frac{\log(2)}{\log(400)}\right)} p_{k}(\alpha_c \log(k)) \ge  0.0398153096.
\end{align*}
If $k \le 399$ and $\l \ge 4$, then using MAGMA we find that
\begin{displaymath}
P_m \ge \sum_{2 \le 2^{\j} \le k} \frac{p_k(2^{\j})}{4}  s_{\neg 2^{\j}}(\l-k) \ge  \frac{0.156250}{4} = 0.0390625,
\end{displaymath}
and thus for all $m$ we have $P_m \ge 0.0390625$ in the symplectic case. 

If $\l=2$, or $3$, then $C=\SL_2(q) \!\times \! \Sp_{d-2}(q)$ and we use Proposition \ref{caseS2part2} to show that the proportion of elements in $C$ of the form $(x,y)$, where $|x|_2=2$ and $|y|_2=1$ is at least $1/16$ when $\l=2$ and $5/128=0.0390625$ when $\l=3$.

In the orthogonal cases, our balanced involution has centralizer
 $\SO^{\epsilon_1}_m(q)\times \SO^{\epsilon_2}_{d-m}(q) \le \SO^{\epsilon}_d(q)$  where $\epsilon=\epsilon_1 \epsilon_2$, and we are assuming that $m \ge 3$. Write $m=2k$ if $m$ is even and $m=2k+1$ if $m$ is odd. If $k \le 399$ and $d \ge 13$ , then Propositions \ref{pro:caseS2parts2j} and \ref{pro:caseO2jp}  and Lemma \ref{lem:notdivisibleby} show that when $k$ is even,
\begin{small}
\begin{align*} 
P_m \ge& \sum_{\substack{2 \le 2^{\j} \le k\\ 2^b \ne k_2}} \frac{p_k(2^{\j})}{4}  s_{\neg 2^{\j}}(\l-k)  +  \left( \frac{p_k(k_2)}{4} - \frac{1}{4k}\right)s_{\neg k_2}(n-k)\\
 &+ \frac{p_k(1)}{4}\frac{(2\l -2k)!}{2^{2\l -2k}((2\l-2k) !)^2} \frac{(2\l-2k-2)}{(2\l-2k-1)} \ge 0.0411987,
\end{align*}
\end{small}
and when $k$ is odd,
\begin{small}
\begin{displaymath}
P_m \! \ge \!\! \sum_{2 \le 2^{\j} \le k}\!\! \frac{p_k(2^{\j})}{4}  s_{\neg 2^{\j}}(\l \! - \!k)  
+ \frac{1}{4}\left( p_k(1)\!-\! \frac{1}{k} \right)\frac{(2\l -2k)!}{2^{2\l -2k}((2\l-2k) !)^2} \frac{(2\l \!- \!2k\!-\!2)}{(2\l \!-\!2k\!-\!1)} \! \ge \!0.0411987. 
\end{displaymath}
\end{small}
If  $(6 \le )d \le 13$ (and $k \le 399$), then $d$ and $m$ are small enough to write down all of the Weyl group elements of $Cl_m(q)$ and $Cl_{d-m}(q)$ and calculate the \tpes of each of their corresponding maximal tori. In Table \ref{tab:smallSO}, we describe the proportions of elements in the Weyl groups of $\SO^{\epsilon}_m(q)$,  for $3 \le m \le 8$, that correspond to maximal tori of \tpe $2^b=2$, $t$, $2t$, and $4t$. It is then straightforward to apply Theorem \ref{thm:maintool} to get sharper lower bounds on $\frac{|Q(2^b)|}{|\GF|}$ than in Proposition \ref{pro:caseO2jp}. An easy case-by-case analysis then shows that if $6\le d \le 13$ (and $m \ge 3$), then $P_m \ge 0.0585$.

It is usually enough to simply find a lower bound on the proportion of elements $(x,y) \in C$ such that $x$ has \tpo $t$ or $2t$, and $y$ has \tpo $1$ or $2$.

\begin{table}[htbp]
\caption{Proportions of Weyl group elements in $W$, when $\GF=\SO_{m}(q)$, corresponding to maximal tori of certain \tpes\hspace{-1mm}.}
\label{tab:smallSO}
\vspace{1mm}
\begin{tabular}{cccccc}
\hline
\vspace*{-4mm}\\
$\GF$ & $2^b$ & \multicolumn{2}{c}{ Proportion of $w \in W$ } & \multicolumn{1}{c}{Lower bound} \\
 &  & \multicolumn{2}{c}{corresponding to maximal } & \multicolumn{1}{c}{for $\frac{|Q(2^b)|}{|\GF|}$} \\
 &  & \multicolumn{2}{c}{tori with \tpe $2^b$}  & \multicolumn{1}{c}{}  \\ \cline{3-4}
 &  & \multicolumn{1}{c}{$q \equiv 1 \imod 4$} & \multicolumn{1}{c}{$q \equiv 3 \imod 4$ \rule{0cm}{3mm}} & \multicolumn{1}{c}{} \\ 
 \vspace*{-4mm}\\
 \hline
$\SO_3(q)$ & $t$\rule{0cm}{3mm}&   $1/2$ &   $1/2$ &   $1/4$ \\
 & {2} &   $1/ 2$ &   $1/ 2$ &   $1/ 4$ \\
$\SO_4^{+}(q)$ & $2t$ &   $1/ 2$ &   $1/ 2$ &   $1/ 4$ \\
 & $t$ &   $1/ 4$ &   $1/ 4$ &   $1/ 8$ \\
 & $2$ &   $1/ 4$ &   $1/ 4$ &   $1/ 8$ \\
$\SO_4^{-}(q)$ & $2t$ & $0$       & $0$       & $0$       \\
 & $t$ &   $1/ 2$ &   $1/ 2$ &   $1/ 4$ \\
 & {2} &   $1/ 2$ &   $1/ 2$ &   $1/ 4$ \\
$\SO_5(q)$ & $2t$ &   $1/ 4$ &   $1/ 4$ &   $1/ 8$ \\
 & $t$ &   $3/ 8$ &   $3/ 8$ &   $3/16$ \\
 & $2$ &   $3/ 8$ &   $3/ 8$ &   $3/16$ \\
$\SO_6^{+}(q)$ & $2t$ &   $1/ 4$ &   $1/ 4$ &   $1/ 8$ \\
 & $t$ &   $1/ 2$ &   $3/ 8$ &   $3/16$ \\
 & $2$ &   $1/ 4$ &   $3/ 8$ &   $1/ 8$ \\
$\SO_6^{-}(q)$ & $2t$ &   $1/ 4$ &   $1/ 4$ &   $1/ 8$ \\
 & $t$ &   $3/ 8$ &   $1/ 2$ &   $3/16$ \\
 & $2$ &   $3/ 8$ &   $1/ 4$ &   $1/ 8$ \\
$\SO_7(q)$ & $2t$ &  $1/ 4$ &   $1/ 4$ &   $1/ 8$ \\
 & $t$ &   $7/16$ &   $7/16$ &   $7/32$ \\
 & $2$ &   $5/16$ &   $5/16$ &   $5/32$ \\
$\SO_8^{+}(q)$ & $4t$ &   $1/ 4$ &   $1/ 4$ &   $1/ 8$ \\
 & $2t$ &   $3/16$ &   $3/16$ &   $3/32$ \\
 & $t$ &  $21/64$ &  $21/64$ &  $21/128$ \\
 & $2$ &  $15/64$ &  $15/64$ &  $15/128$ \\
$\SO_8^{-}(q)$ & $4t$ & $0$       & $0$       & $0$       \\
 & $2t$ &   $1/ 4$ &   $1/ 4$ &   $1/ 8$ \\
 & $t$ &   $7/16$ &   $7/16$ &   $7/32$ \\
 & $2$ &   $5/16$ &   $5/16$ &   $5/32$ \\
 \vspace*{-4mm}\\
  \hline
\end{tabular}
\end{table}

 If $k\ge 400$, then
\begin{align*}
P_m& \! \ge \! \sum_{[c,2c) \in \I} s_{\neg \alpha_c \log(k)}(\l-k) p_{k}(\alpha_c \log(k))/4  - \frac{1}{4(400)} \ge 0.039190,
\end{align*}
hence $P_m \ge  0.039190$ for all $m$ in the orthogonal cases.
%

Now if we suppose that $ d/2 \le m \le 2d/3$, then the same calculation yields that if $m\ge 400$, then
\begin{align*}
s_{\neg \alpha_c \log(m)}(n-m) &\ge c(\alpha_c \log(m))(1-1/(n-m))(n-m)^{- \frac{1}{\alpha_c \log(m)}}\\
&\ge c(\alpha_c \log(m))(1-2/m)m^{- \frac{1}{\alpha_c \log(m)}} \\
& \ge  c(\alpha_c \log(m))(398/400)e^{-\frac{1}{\alpha_c}}.
\end{align*}
So we can use, in the linear case for example, the inequality
\begin{small}
\begin{displaymath}
P_m \ge \sum_{[c,2c) \in \I} \frac{p_m(\alpha_c \log(m))}{2}  s_{\neg \alpha_c \log(m)}(n-m) \ge 0.175366870998411
/2 =0.08768343549921.
\end{displaymath}
\end{small}
We can verify in MAGMA that this lower bound also holds when $m \le 400$. The other cases are similar.
\end{proof}

\section{Analysis of Yal\c{c}inkaya's algorithm} \label{sec:sukru}
Next we consider an involution in a classical group with fixed space of arbitrary dimension (at least $2$), and prove the following theorem, which then allows us to prove Corollary \ref{cor:sukru}. Again we note that for an integer $m$ and $a$ as in \eqref{e:a}, the interval $\left[ \frac{1}{2}a\log(m), a\log(m)\right)$ contains a unique power of $2$, which we denote by $p_{0,m}$.
\begin{thm} \label{thm:sukru}
Let $\GF = \GL^{\epsilon}_{n}(q)$, $\Sp_{2n}(q)$, $\SO_{2n+1}(q)$, or $\SO^{\epsilon}_{2n}(q)$, where $q$ is odd, and $d$ is the dimension of the natural module for $G$.
Consider an involution in $G$ with centralizer $C = Cl_m(q) \times Cl_{d-m}(q)$.  Suppose that $d/2,m \ge 2$, and $Cl_m(q) \ne \SO_2^{\epsilon}(q)$. Let $P^{\GF}_m$ and $a$ be as in \eqref{e:PMG} and \eqref{e:a}. If $\GF = \GL^{\epsilon}_n(q)$, $m\ge 16$, and $p_{0}=p_{0.m}$ then
\begin{align}
P_m^{\GF}\ge &   \frac{1}{1000}\left (15s_{\neg 8 p_0}(n\!-\!m)\!+\! 30s_{\neg 4 p_0}(n\!-\!m) \! + \!76s_{\neg 2p_0}(n\!-\!m)\!+\!117s_{\neg p_0}(n\!-\!m)\right).\notag 
\end{align}
If $\GF = \GL^{\epsilon}_n(q)$,  $m\le 15$ and $p_0=p_{0,m}$, then Table \ref{tab:PmGL} describes lower bounds for $P_m^{\GF}$.
If $\GF = \Sp_{2\l}(q)$, $\SO_{2n+1}(q)$, or $\SO^{\epsilon}_{2n}(q)$ and if $k:= \left \lfloor \frac{m}{2}\right\rfloor \ge 17$ and $p_0 = p_{0,k}$, then
\begin{align}
P_m^{\GF}\ge& \frac{1}{500}\left(2s_{\neg 8 p_0}(\l-k)+ 7s_{\neg 4 p_0}(\l-k) +19s_{\neg 2p_0}(\l-k)+29s_{\neg p_0}(\l-k)\right ). \notag
\end{align}

If $\GF = \Sp_{2\l}(q)$, $\SO_{2n+1}(q)$, or $\SO^{\epsilon}_{2n}(q)$, $k \le 16$ and $p_0=p_{0,k}$, then Tables \ref{tab:PmSpSO} and \ref{tab:PmSO} describe lower bounds for $P_m^{\GF}$.
\end{thm}
\begin{table}
\caption{Lower bounds for $P_m^\GF$ when $G=\GL^{\epsilon}_{n}(q)$  and $m\le 15$.} \label{tab:PmGL}
\vspace{1mm}
\begin{tabular}{cc}
  \hline
  \vspace*{-4mm}\\
  $m$ & $P_m^G$ \rule{0cm}{3mm}\\
  \vspace*{-4mm}\\
    \hline
  \vspace*{-4mm}\\
  $2$, $3$ & $\frac{ s_{\neg 2}(n-m)}{4}$ \rule{0cm}{3.5mm} \\
  $4$, $5$ & $\frac{ s_{\neg 4}(n-m)}{8} + \frac{3 s_{\neg 2}(n-m)}{16}$ \\
  $6$, $7$ & $\frac{ s_{\neg 4}(n-m)}{8} + \frac{7 s_{\neg 2}(n-m)}{32}$ \\
  $8$, $9$, $10$, $11$ &  $\frac{ s_{\neg 8}(n-m)}{16}+  \frac{7 s_{\neg 4}(n-m)}{64} + \frac{49 s_{\neg 2}(n-m)}{256}$ \\
 $12$, $13$, $14$, $15$  & $\frac{ s_{\neg 8}(n-m)}{16}+  \frac{35 s_{\neg 4}(n-m)}{256} + \frac{385 s_{\neg 2}(n-m)}{2048}$\rule[-4pt]{0cm}{3mm} \\
 \vspace*{-3.5mm} \\
 \hline
\end{tabular}
\end{table}

\begin{table}
\caption{Lower bounds for $P_m^G$ when $G=\Sp_{2\l}(q)$ and $k= \frac{m}{2} \le 16$.} \label{tab:PmSpSO}
\vspace{1mm}
\begin{tabular}{cc}
  \hline
  \vspace*{-4mm}\\
  $k$ & $P_m^G$ \rule{0cm}{3mm}\\
  \vspace*{-4mm}\\
  \hline
  $1$      & $\frac{25}{116 \sqrt{\pi (n-1)}}$ (see Remark \ref{rem:SL2})\rule{0cm}{3.5mm}\\
  $2$, $3$ & $\frac{ s_{\neg 2}(\l-k)}{8} + \frac{25 }{232 \sqrt{\pi (\l-k)}}$ \\
  $4$, $5$ & $\frac{ s_{\neg 4}(\l-k)}{16} + \frac{3 s_{\neg 2}(\l-k)}{32}+\frac{ 75 }{928 \sqrt{\pi (\l-k)}}$ \\
  $6$, $7$ & $\frac{ s_{\neg 4}(\l-k)}{16} + \frac{7 s_{\neg 2}(\l-k)}{64}+\frac{125 }{1856 \sqrt{\pi (\l-k)}}$ \\
  $8$, $9$, $10$, $11$ &  $\frac{ s_{\neg 8}(\l-k)}{32}+  \frac{7 s_{\neg 4}(\l-k)}{128} + \frac{49 s_{\neg 2}(\l-k)}{512} + \frac{1575 }{29696 \sqrt{\pi (\l-k)}}$ \\
   $12$, $13$, $14$, $15$, $16$& $\frac{ s_{\neg 8}(\l-k)}{32}+  \frac{35 s_{\neg 4}(\l-k)}{512} + \frac{385 s_{\neg 2}(\l-k)}{4096}  +\frac{10725 }{237568 \sqrt{\pi (\l-k)}}$\rule[-9pt]{0cm}{9mm} \\
   \vspace*{-4mm}\\
  \hline
\end{tabular}
\end{table}

\begin{table}
\caption{Lower bounds for $P_m^G$ when $G=\SO_{d}^{\epsilon}(q)$, $m \le 33$. Here $\chi(\l-k)=1$ if $d-m$ is odd, and $\chi(\l-k) \ge \frac{2\l-2k-2}{2\l-2k-1}$ otherwise.} \label{tab:PmSO}
\vspace{1mm}
\begin{tabular}{cc}
  \hline
  \vspace*{-4mm}\\
  $k$ & $P_m^G$ \rule{0cm}{3mm}\\
  \vspace*{-4mm}\\
  \hline
  \vspace*{-4mm}\\
  $1$ and $m=3$ & $\frac{25}{116 \sqrt{\pi(\l-2)}}$ \rule{0cm}{3.5mm}\\
  $2$ & $\frac{25 \chi(\l-k)}{232 \sqrt{\pi (\l-k)}}$ \\
  $3$ & $\frac{ s_{\neg 2}(\l-k)}{8} + \frac{25 \chi(\l-k)}{232 \sqrt{\pi (\l-k)}}$ \\
$4$ &  $\frac{3 s_{\neg 2}(\l-k)}{32}+\frac{ 75 \chi(\l-k)}{928 \sqrt{\pi (\l-k)}}$ \\
  $5$ & $\frac{ s_{\neg 4}(\l-k)}{16} + \frac{3s_{\neg 2}(\l-k)}{32}+\frac{75 \chi(\l-k)}{928 \sqrt{\pi(\l-k)}}$ \\
  $6$, $7$ & $\frac{s_{\neg 4}(\l-k)}{16} + \frac{13 s_{\neg 2}(\l-k)}{192}+\frac{125 \chi(\l-k)}{1856 \sqrt{\pi (\l-k)}}$ \\
 $8$&   $\frac{7s_{\neg 4}(\l-k)}{128} + \frac{49s_{\neg 2}(\l-k)}{512} + \frac{1575 \chi(\l-k)}{29696 \sqrt{\pi (\l-k)}}$ \\
$9$, $10$, $11$&  $\frac{s_{\neg 8}(\l-k)}{32}+  \frac{7 s_{\neg 4}(\l-k)}{128} + \frac{397 s_{\neg 2}(\l-k)}{5120} + \frac{1575 \chi(\l-k)}{29696 \sqrt{\pi (\l-k)}}$ \\
  $12$, $13$, $14$, $15$, $16$  & $\frac{ s_{\neg 8}(\l-k)}{32}+  \frac{35 s_{\neg 4}(\l-k)}{512} + \frac{385 s_{\neg 2}(\l-k)}{4096}  +\frac{10725 \chi(\l-k)}{237568 \sqrt{\pi (\l-k)}}$\rule[-9pt]{0cm}{9mm} \\
  \vspace*{-4mm}\\
  \hline
\end{tabular}
\end{table}

 \begin{remark}\label{rem:SL2}
 \emph{
If $k=1$ and $\GF=\Sp_{2n}(q)$, then $C=\SL_2(q) \times \Sp_{2n-2}(q)$, and Theorem \ref{thm:sukru} gives a lower bound for the proportion of elements $(x,y) \in C$ that power up to $(w,1)$, where $w$ is the central involution in $\SL_2(q)$. By taking $x$ of \tpo $t$ and $y$ of odd order, we can show that the proportion of elements $(x,y) \in C$ that power up to $(w,1)$, where $w$ has order $4$ in $\SL_2(q)$, is at least $\frac{\Gamma(n-3/4)}{4\Gamma(1/4)(n-1)!} \ge \frac{1}{16n^{3/4}}$.}
 \end{remark}
\begin{proof}
The proof is similar to the proof of Proposition \ref{pro:balanced}. In the linear and unitary cases, we make use of the inequality
\begin{align} \label{eqn:last}
P_m &\ge \sum_{2 \le 2^{\j} \le m} \frac{|Q_{m}(2^{\j}t)|}{|\GF|}  s_{\neg 2^{\j}}(n-m), \end{align}
where $Q_m(2^{\j}t)$ is the set of elements in $Cl_m(q)$ of \tpo $2^{\j}t$.
The cases where $m\in \{2, \ldots, 15\}$ follow since $\frac{|Q_{m}(2^{\j}t)|}{|\GF|}\ge \frac{p_m(2^{\j})}{2}$ and we calculate $p_m(2^{\j})$ explicitly.
 The case $m \ge 16$ follows in the same way except that we estimate $p_m(2^{\j})$ using Theorem \ref{thm:constantSn}, where we only take the terms in the sum corresponding to $2^{\j} = p_0,2 p_0, 4 p_0, 8 p_0$ and $p_0 \in [ \frac{1}{2}a \log(m), a \log(m)]$.

In the symplectic case we use a similar inequality to the one in the linear case, but with an extra term, representing the proportion of elements with \tpo $t$ in $Cl_m(q)$ and \tpo $1$ or $2$ in $Cl_{d-m}(q)$. We have
\begin{displaymath}
P_m \ge \sum_{2 \le 2^{\j} \le k} \frac{|Q_{m}(2^{\j}t)|}{|\GF|}  s_{\neg 2^{\j}}(\l-k) +  \frac{25}{29\sqrt{\pi(\l -k)}}\frac{s_{\neg 2}(k)}{4}, \end{displaymath}
which is justified using Propositions \ref{caseS2parts} and \ref{caseS2part2}. We use the lower bounds on $\frac{|Q_{m}(2^{\j}t)|}{|\GF|}$ of Propositions \ref{caseS2parts} and \ref{pro:caseS2parts2j}, which are given in terms of $p_k(2^{\j})$, and we calculate $p_k(2^{\j})$ explicitly for $k \le 16$. For $k \ge 17$, we use the constant lower bounds on $\frac{|Q_{m}(2^{\j}t)|}{|\GF|}$ given in Table \ref{tab:forpropenot2e}. The orthogonal case is similar.
\end{proof}
 \begin{proof}[Proof of Corollary \ref{cor:sukru}]
 Corollary \ref{cor:sukru} now follows since Theorem \ref{thm:sukru} shows that there exists a positive constant $K_1$ such that $P_m \ge K_1 s_{\neg 2^\a}(n-m)$ or $P_m \ge K_1 s_{\neg 2^\a}(n-k)$ for some positive integer $\a$ depending on $m$ or $k$. But $s_{\neg 2^\a}(n-m)\ge s_{\neg 2}(n-m)$ and $s_{\neg 2^\a}(n-k)\ge s_{\neg 2}(n-k)$. Moreover, Theorem \ref{thm:constantSn}(c) shows that there exists a positive constant $K_2$ such that $s_{\neg 2}(n-m) \ge K_2/ \sqrt{n}$ whenever $n-m \ge 2$ and $s_{\neg 2}(n-k) \ge K_2/ \sqrt{n}$ whenever $n-k \ge 2$.
 \end{proof}

\subsection*{Acknowledgements}
The authors would like to thank Tomasz Popiel for several preliminary discussions. They would like to thank Pablo Spiga for pointing out the significance of the Stirling numbers in counting elements of the symmetric group, which led to Theorem \ref{thm:oddandtwiceodd}. They would also like to thank \c{S}\"{u}kr\"{u} Yal\c{c}inkaya for several discussions regarding his algorithm, which led to the results in Section \ref{sec:sukru}.

\bibliographystyle{amsalpha}
\bibliography{/Users/simon_guest/Dropbox/mbibliography2}
\end{document}